\documentclass[review,hidelinks,onefignum,onetabnum]{siamart220329}
\usepackage{pgfplots}
\usepackage{caption}
\usepackage{lipsum}
\usepackage{amsfonts}
\usepackage{graphicx}
\usepackage{epstopdf}
\usepackage{algorithmic}
\ifpdf
  \DeclareGraphicsExtensions{.eps,.pdf,.png,.jpg}
\else
  \DeclareGraphicsExtensions{.eps}
\fi

\newsiamremark{remark}{Remark}
\newsiamremark{hypothesis}{Hypothesis}
\crefname{hypothesis}{Hypothesis}{Hypotheses}
\newsiamthm{claim}{Claim}

\headers{An Example Article}{D. Doe, P. T. Frank, and J. E. Smith}

\title{An Example Article\thanks{Submitted to the editors DATE.
\funding{This work was funded by the Fog Research Institute under contract no.~FRI-454.}}}

\author{Dianne Doe\thanks{Imagination Corp., Chicago, IL 
  (\email{ddoe@imag.com}, \url{http://www.imag.com/\string~ddoe/}).}
\and Paul T. Frank\thanks{Department of Applied Mathematics, Fictional University, Boise, ID 
  (\email{ptfrank@fictional.edu}, \email{jesmith@fictional.edu}).}
\and Jane E. Smith\footnotemark[3]}

\usepackage{amsopn}

\nolinenumbers
\ifpdf
\hypersetup{
  pdftitle={GENERALIZED PATCH DYNAMICS SCHEME IN EQUATION-FREE MULTISCALE MODELLING},
  pdfauthor={T. K. KARMAKAR, and D. C. DALAL}
}
\fi

\usepackage{braket,amsfonts}

\usepackage{array}

\usepackage[caption=false]{subfig}
\usepackage{pgfplots}

\crefname{hypothesis}{Hypothesis}{Hypotheses}

\usepackage{algorithmic}

\usepackage{graphicx,epstopdf}
\usepackage{array,multirow}

\usepackage{amsopn}

\usepackage{xspace}
\usepackage{bold-extra}
\usepackage[most]{tcolorbox}

\colorlet{texcscolor}{blue!50!black}
\colorlet{texemcolor}{red!70!black}
\colorlet{texpreamble}{red!70!black}
\colorlet{codebackground}{black!25!white!25}

\usepackage{graphicx,epstopdf} % <- Preamble
\usepackage[caption=false]{subfig} % <- Preamble
\usepackage{subfig}

\begin{tcbverbatimwrite}{tmp_\jobname_header.tex}
	\title{GENERALIZED PATCH DYNAMICS SCHEME IN EQUATION-FREE MULTISCALE MODELLING\thanks{Submitted to the editors January 1, 2024.
			\funding{Current research did not receive any fund from any external source.}}}
	
	\author{T. K. Karmakar\thanks{Department of Mathematics, Indian Institute of Technology Guwahati, Guwahati, Assam 781039, India (\email{tanay.kumar@iitg.ac.in}; \email{durga@iitg.ac.in}).}
		\and D. C. Dalal\footnotemark[2]}
	
	\headers{GENERALIZED PATCH DYNAMICS SCHEME}{T. K. Karmakar, and D. C. Dalal}
\end{tcbverbatimwrite}
	\title{GENERALIZED PATCH DYNAMICS SCHEME IN EQUATION-FREE MULTISCALE MODELLING\thanks{Submitted to the editors January 1, 2024.
			\funding{Current research did not receive any fund from any external source.}}}
	
	\author{T. K. Karmakar\thanks{Department of Mathematics, Indian Institute of Technology Guwahati, Guwahati, Assam 781039, India (\email{tanay.kumar@iitg.ac.in}; \email{durga@iitg.ac.in}).}
		\and D. C. Dalal\footnotemark[2]}
	
	\headers{GENERALIZED PATCH DYNAMICS SCHEME}{T. K. Karmakar, and D. C. Dalal}

\begin{document}

\maketitle

% REQUIRED
\begin{abstract}
	There is a class of problems that exhibit smooth behavior on macroscopic scales, where only a microscopic evolution law is known. Patch dynamics scheme of ‘equation-free multiscale modelling’ is one of the techniques, which aims to extract the macroscopic information using such known time-dependent microscopic model simulation in patches (which is a fraction of the space-time domain)  that  reduces  the  computational  complexity. Here, extrapolation time step has an important role to reduce the error at macroscopic level. In this study, a generalized patch dynamics (GPD) scheme is proposed by distributing the gap-tooth timesteppers (GTTs) within each long (macroscopic) time step. This distribution is done in two ways, namely, GPD schemes of type-I and type-II. The proposed GPD scheme is based on three different time scales namely, micro, meso and macro to predict the system level behaviours. The GPD scheme of both types are capable of providing better accuracy with less computation time compared to the usual patch dynamics (UPD) scheme. The physical behaviours of the problems can be more appropriately addressed by the GPD scheme as one may use a non-uniform (variable) distribution of gap-tooth timesteppers (GTTs), as well as the extrapolation times based on the physics of the problem. Where the UPD scheme fails to converge for a long extrapolation time, both types of GPD schemes can be successfully applied. The whole method has been analyzed successfully for the one-dimensional reaction-diffusion problem.
\end{abstract}

% REQUIRED

% REQUIRED

\section{Introduction}
The behaviour of a system on a large domain can be understood through a fine-level simulation over the entire physical domain, which is a classical numerical concept. But computation time and memory constraints make such simulation futile. A balanced model should be computationally feasible without losing much information; such a model can be obtained through multiscale modelling. In science and engineering, the governing microscopic mechanisms are usually known, but translation of microscopic description to system level macroscopic description is rarely available. For a class of multiscale problems, separation of scales appear between the microscopic model and the macroscopic model. Main challenge of multiscale modelling is to form scale bridging to transform and transfer the information between the different scales. In equation-free multiscale modelling, the available microscopic model must be closed in the whole domain, but macroscopic model need not be closed. Equation-free computation aims to extract the macroscopic information using the given microscopic model simulation in a fraction of the space-time domain that reduces the computational complexity. This paper sheds light on the generalization of the patch dynamics scheme in equation-free multiscale modelling.

Equation-free framework \cite{2003_kevrekidis_equation} is built on coarse time-stepper, which is a small macroscopic time mapping from coarse variables to coarse variables. Coarse time-stepper consists of three steps: (i) lifting, i.e., transforming the macroscopic initial data into an appropriate initial condition for the microscopic model; (ii) evolution, i.e., evolving the system for short macroscopic time with the help of microscopic simulator; (iii) restriction, i.e., use of restriction operator on the microscopic variables to get back the macroscopic variables. The above three steps finally produce coarse time-stepper solutions. A coarse time-stepper can be used with coarse projective integration (CPI) method for large time integration \cite{2003_gear_projective}. Time-stepper based bifurcation analysis can be used to find the coarse solution at the infinite time for the unavailable macroscopic equation \cite{2002_gear_coarse,2002_makeev_coarse,2002_makeev_coarse1}. This approach has several applications \cite{2003_hummer_coarse,2003_siettos_coarse}.

For a system having spatial and temporal variables, Gear et al. \cite{2003_gear_gap,2003_kevrekidis_equation} proposed the gap-tooth scheme (GTS) for the system with one space dimension and subsequently it was extended to higher space dimensions. The GTS is the generalization of the coarse time-stepper, where space variables are involved. GTS bridges ``small space, short time" to ``large space, short time" evolution. The combination of GTS and CPI is called as patch dynamics
scheme \cite{2003_kevrekidis_equation,2009_samaey_equation}. It provides an equation-free bridging between “small space, short time” to “large space, long
time” evolution. The patch dynamics scheme used so far in the literature  will be termed as usual patch dynamics (UPD) scheme now on. 

%In order to improve the efficiency of long time-step integration, one can combine GTS and CPI. In order to combine this a few number of GTTs each of short time $\tau$ are used and subsequently the time derivative of the unavailable macroscopic equation is estimated. This estimate is used to perform a long time integration of size $\Delta t\gg\tau$ to reach the next coarse level. The combination of GTS and CPI is called as patch dynamics scheme  \cite{2003_kevrekidis_equation,2009_samaey_equation}. \textbf{To distinguish it from the proposed scheme, let's name it the usual patch dynamics (UPD) scheme.} It provides an equation-free bridging between ``small space, short time" to ``large space, long time" evolution.      

When more information about the unavailable macroscopic problem is known, it helps to modify the general time-stepper based procedure appropriately. The coupling of the microscopic simulator with macroscopic behaviour needs information about the nature of the unavailable coarse equation. Li et al. \cite{2007_li_deciding} have proposed ``The Baby-Bathwater Scheme" to decide the nature of the coarse equation through microscopic simulation. 

In some articles \cite{2003_kevrekidis_equation,2005_samaey_gap,2006_samaey_patch,2020_arbabi_linking,2015_liu_acceleration,2009_samaey_equation}, forward Euler's method is used to extrapolate the macroscopic solution for large time-step as outer integrator in projective integration. Later, new projective version of second-order accurate Runge-Kutta and Adams-Bashforth methods were introduced by Lee et al. \cite{2007_lee_second} as an outer integrator to solve stiff differential systems. Maclean et al. \cite{2015_maclean_convergence} presented a convergence analysis for higher order projective integration method, stability condition and error bounds for slow variables for a class of multiscale systems. 

The fundamental idea of patch dynamics scheme is to replace the expensive simulation across a large domain with sparse computation across small, well-separated patches within that domain. Patch dynamics depends on several parameters, like width of the patch ($h$), macroscopic and microscopic spatial domain discretization ($\Delta x$ \& $\delta x$), short time step ($\tau$) for gap-tooth timestepper, number of gap-tooth timesteppers ($k$), long macroscopic time step size ($\Delta t$), microscopic time discretization ($\delta t$) etc. For best results, it is crucial to choose the appropriate patch parameters. The patch scheme is accurate up to order $\mathcal{O}((\Delta x)^{2\Gamma})$, where $\Delta x$ is the distance between patches and $\Gamma$ is the order of coupling between patches \cite{2007_roberts_general,2014_roberts_dynamical,2016_cao_multiscale}, when the microscale is homogeneous or heterogeneous \cite{2017_bunder_good}. Equation-Free Toolbox \cite{2021_maclean_toolbox} gives a practical introduction of patch dynamics scheme. By restricting the transfer of data between processors, Bunder et al. \cite{2016_bunder_accuracy} modified the patch dynamics macroscale modeling scheme for massive parallelization and exascale computing. In this modified scheme, the concept of meso scale becomes evident.

Problems involving highly variable microscale diffusion with significant roughness in the microscale structure can also be solved using the equation-free scheme. Bunder et al. \cite{2011_bunder_patch,2017_bunder_good} proposed a patch dynamics scheme for 1D heterogeneous diffusion equation at the microscopic level and it was shown that the accurate patch coupling conditions take care of the underlying microscale structure, and the results show this method can be applied to a range of microscale heterogeneity issues.

To perform time integration in the multiscale problems, both micro and macro time steps play crucial roles. The main objective of this paper is to propose a general framework for designing multiscale, multi-physics methods. This general framework should have distinct advantages over other existing general stategies,  such as systematic up-scaling \cite{1977_brandt_multi,2002_brandt_multiscale}, the heterogeneous multiscale method (HMM) \cite{2003_engquist_heterognous,2003_vanden_numerical,2007_engquist_heterogeneous} and the ``equation-free" approach \cite{2003_kevrekidis_equation,2006_samaey_patch,2020_arbabi_linking}. 

One of the improved versions of the HMM is the seamless heterogeneous multiscale method (SHMM), first introduced by Fatkullin et al. \cite{2004_fatkullin_computational} and further improved by W. E. \cite{2009_weinan_general}, where reinitialization of microscale simulation is not required at each macro time step or each macro iteration step.

In all preceding studies, attention was solely directed towards two distinct time-scales, namely micro and macro time scales, to capture the effective behaviour over a long time. In this paper, a new time-scale is proposed, the meso time scale, which lies between micro and macro time-scales. The macroscale solver uses mesoscale time step finer than the macro time scale to find the system level behaviours. This new multiscale integrator can be applied to handle stiff and highly oscillatory problems more efficiently. In order for the microscale system to relax and influence the macroscale system, the macroscale solver uses a mesoscale time step, much finer than the one in HMM, known as mesoscale HMM (MSHMM) \cite{2009_weinan_general,2013_lee_variable}.

Without identifying the slow and fast variables in advance, a similar approach can also be applied to solve the system. This idea was initially proposed by Vanden-Eijnden \cite{2007_engquist_heterogeneous} and subsequently modified by Tao et al. \cite{2010_tao_nonintrusive}, known as `Flow Averaging Integrators' (FLAVORS). In FLAVORS, the stiff part is turned on over microscopic time step and off over mesoscopic time steps. 

To handle the transient and amplified oscillations that are generated in UPD/HMM/SHMM schemes, the non-uniform GPD schemes of both types are a better choice to model the multiscale problems and to preserve the computational complexity of the method. To reduce the error in an efficient way, the GPD scheme uses variable mesoscopic time steps (i.e., variable number of GTTs and variable extrapolation times) with respect to the physics of the problem. However, consecutively, the computational complexity increases exponentially as the number of different scales increases \cite{2015_lee_fast}. This kind of behaviours in the solution can also be observed in the GPD scheme. The stiffness parameter ($\epsilon$) is a sensative parameter in highly stiff and oscillatory problems. W. E. et al. \cite{2009_weinan_general} and Maclean \cite{2015_maclean_note} proposed the `Boosting Algorithm (BA)' in HMM. By increasing the value of $\epsilon$, the algorithm can be boosted, and the system can be made less stiff, resulting in significant computational advantages.

The idea of variable time step is also seen in `variable step size heterogeneous multiscale method (VSHMM)' proposed by Lee et al. \cite{2013_lee_variable}, but the method was only restriced to ODEs. VSHMM was later extended to handle several temporal scale problems without iteratively applying two-scale integrators \cite{2015_lee_fast}. In \cite{2013_lee_variable}, it employs very fine mesoscopic time steps at the beginning and the end of each macro time step to overcome problems such as the delayed relaxation of the fast variables in dissipative problems and amplified oscillations in	highly oscillatory problems. Meanwhile, it uses coarse mesoscopic time steps at all other times to save computational complexity. Once the system has evolved to the next macro time step, the same process is iterated. Such features are not available in the UPD scheme, but are present in the GPD scheme. These uniform or variable mesoscopic step sizes contribute to making the model more realistic than the UPD scheme.

The concept behind the GPD scheme is a general framework which allows for the modification of the systematic up-scaling and the HMM techniques using a unified strategy. By maintaining the same stucture as the GPD scheme (figure \ref{fig:3}, \ref{fig:4}), the systematic up-scaling and the HMM techniques can be applied on it based on their own macro to micro operator, micro to macro operator, inner \& outer integrators, and bridging techniques. Existing methods HMM \cite{2003_engquist_heterognous,2003_vanden_numerical,2007_engquist_heterogeneous}, FLAVORS \cite{2010_tao_nonintrusive}, VSHMM \cite{2013_lee_variable}, BA \cite{2009_weinan_general,2015_maclean_note}, UPD scheme \cite{2003_kevrekidis_equation,2006_samaey_patch,2020_arbabi_linking} are the particular cases of the mother scheme i.e., the GPD scheme.

%\textbf{$\underline{Boosting Algorithm:}$ 
%By putting $l=1$ \& $k=0$ in the non-uniform GPD scheme of type-I with `Boosting Algorithm (BA)', one can obtain the `BA in HMM' as proposed by John Maclean \cite{2015_maclean_note}. }

In this paper, a generalized version of the patch dynamics (GPD) scheme is proposed by distributing the GTTs within each long macroscopic time step $\Delta t$ which separates the time scale into three different scales, which are microscale, mesoscale and macroscale. This distribution can be done in two ways, which are named as GPD schemes of type-I and type-II. The proposed distribution of the GTTs decreases the extrapolation time length and improves the solution. Pre-existing schemes can be obtained from the GPD scheme by imposing some particular values. These GPD schemes are capable of producing a better accuracy in a lesser computataional time compared to the UPD scheme. GPD schemes are more suitable for efficiently capturing the physics of the problems compared to the existing methods, as it incorporate additional features that contribute to accurately representing the physical behavior. One of the most important features is that one has the independence of choosing uniform or non-uniform distribution of the GTTs, as well as the extrapolation times with respect to corresponding physics of the problem. After a certain limit of long extrapolation time, the UPD schemes mostly diverge. In most of such cases, the GPD schemes can easily overcome this kind of drawbacks of UPD scheme.

\section{Equation-Free Multiscale Modelling}
\label{sec:alg}
Let the coarse evolution equation in the one-dimensional spatial domain [0,1] be
\begin{equation}\label{eqn:1}
\frac{\partial U}{\partial t}=F(U,\partial_xU,\partial_x^2U,...,\partial_x^dU,x,t),
\end{equation}
in which $U(x,t)$ denotes the macroscopic state, $``d"$ denotes the order of the partial differential equation \eqref{eqn:1}, $\partial_x^kU$ denotes the $k^{th}$ spatial derivative $(k=1,2,...,d )$.

Consider the unsteady microscopic problem in one-dimensional domain [0,1] as

\begin{equation}\label{eqn:2}
\frac{\partial u}{\partial t}=f(u,\partial_xu,\partial_x^2u,...,\partial_x^du,x,t),
\end{equation}
which is in closed form.

The space domain $[0,1]$ is discretized using an equidistant macroscopic mesh $\Pi(\Delta x)$ := $\{x_0=0 < x_1< ... <  x_i= x_0+i\Delta x< ... < x_N = 1\}$ so that the closed approximation is obtained. A method-of-line space discretization of the given equation \eqref{eqn:1} is defined as,

\begin{equation}\label{eqn:3}
\frac{\partial U_j(t)}{\partial t}=F\left(U_j(t), D^1(U_j(t)),...,D^d(U_j(t)),t\right), \hspace{1cm} j=1,...,N-1,
\end{equation}
where, $U_j(t)\approx U(x_j,t)$ and $D^k(U_j(t))$ denotes finite difference approximation for the $k^{th}$ spatial derivative, $\forall k=1,2,...,d$. 

\subsection{Boundary Conditions}

It is very important to choose the suitable patch boundary conditions based on the behaviours of the system in small patches. The macroscopic field $U$ locally assimilates a polynomial between teeth, i.e.,

\begin{equation}\label{eqn:5}
u(x,t_n+m\tau)\approx p_j^\eta(x,t_n+m\tau), \hspace{0.5cm} x\in\left[x_j-\frac{h}{2},x_j+\frac{h}{2}\right],
\end{equation}
where, $m=0, 1,..., k$ and $p_j^\eta$ denotes a polynomial of even degree $\eta$ within $j^{th}$ tooth, $\forall$ $j=1,...,N-1$. The coefficients of the approximating polynomial can be determined by using the concept of box averages in the $j^{th}$ box, in $\frac{\eta}{2}$ boxes to the left and same number to the right of the $j^{th}$ box:

\begin{equation}\label{eqn:6}
\frac{1}{h}\int_{x_{i+j}-\frac{h}{2}}^{x_{i+j}+\frac{h}{2}}p_j^\eta(\xi,t_n+m\tau)d\xi={U}_{i+j}^{m,n}, \hspace{0.5cm} i=-\frac{\eta}{2},...,\frac{\eta}{2}.
\end{equation}

It is easy to verify

\begin{equation}
\mathcal{S}_h(p_j^\eta)(x,t_n+m\tau)= \sum_{i=-\frac{\eta}{2}}^{\frac{\eta}{2}}U_{i+j}^{m,n}L_{i,j}^\eta(x), \hspace{0.25cm} L_{i,j}^\eta(x)=\prod_{s=-\frac{\eta}{2}}^{\frac{\eta}{2}}\frac{(x-x_{j+s})}{(x_{i+j}-x_{j+s})} \hspace{0.2cm}\text{with}\hspace{0.2cm} (s\neq i)
\end{equation}
where $L_{i,j}^\eta$ is a Lagrange polynomial of degree $\eta$ and $\mathcal{S}_h$ is an averaging operator for micro variable $u$ in small patches of size $h$. $U^{m,n}$ is the approximate value of the coarse variable $U$ at time $t_n+m\tau$.

In the GTS, higher order accuracy can be achieved on macroscale for any tooth boundary condition in the microscopic simulator corresponding to the physics of the problem \cite{2007_roberts_general}. For Neumann boundary conditions, the slope of the polynomial for simulating a detailed solution within the $j^{th}$ tooth, during the entire duration of the next time step, 

\begin{equation}\label{eqn:7}
\frac{\partial u}{\partial x}\bigg|_{x_j\pm \frac{h}{2}}= \frac{\partial p_j^\eta}{\partial x}\bigg|_{x_j\pm \frac{h}{2}}=s_j^{\pm}, \hspace{0.5cm} t \in \left[t_n+m\tau, t_n+(m+1)\tau\right].
\end{equation}

\subsection{Initial condition}

For time integration, an initial condition $u_j(x,t_n+m\tau)$ must be needed in each box $\left[x_j-\frac{h}{2},x_j+\frac{h}{2}\right]$, at time $t_n+m\tau$. A standard choice for lifting operator which is a polynomial expansion within the $j^{th}$ tooth can be written as

\begin{equation}\label{eqn:8}
u_j(x,t_n+m\tau)=\sum_{k=0}^{d}D_j^k(U^{m,n})\frac{(x-x_j)^k}{k!}, \hspace{0.5cm} x\in \left[x_j-\frac{h}{2},x_j+\frac{h}{2}\right], 
\end{equation}
where $m\in \{0, 1,..., k\}$ and `$d$' is the order of the macro PDE. $D_j^k(U^{m,n})$ is the finite difference approximation for the $k^{th}$ spatial derivative at $(x_j,t_n+m\tau)$. Constant term, $D_j^0(U^{m,n})$ is obtained from the box averaging property,

\begin{equation}\label{eqn:9}
\frac{1}{h}\int_{x_{j}-\frac{h}{2}}^{x_{j}+\frac{h}{2}}u_j(\xi,t_n+m\tau)d\xi=U_{j}^{m,n},
\end{equation} 
as $D_j^0(U^{m,n})=U_j^{m,n}-\frac{1}{24}h^2D_j^2(U^{m,n})$.

\section{The algorithm of usual patch dynamics}
\label{sec:alg_UPD}
Patch Dynamics is a combination of GTS and CPI. The complete algorithm to progress large time step from one coarse time level $t_n$ to the next coarse time level $t_n+\Delta t$ is given below:

\textbf{Gap-tooth Scheme:}

\textbf{(i) Boundary conditions:} The values needed in the patch boundary conditions is computed from the coarse field using \eqref{eqn:5}. Here, Neumann boundary conditions \eqref{eqn:7} are considered for each and every patch. 

\textbf{(ii) Lifting:} At time $t_n$, the microscopic initial condition $u_j(x,t_n)$ is constructed inside the patch using \eqref{eqn:8}. 

\textbf{(iii) Evolution:} The micro model \eqref{eqn:2} is computed inside the box $\left[x_j-\frac{h}{2},x_j+\frac{h}{2}\right]$ with the boundary conditions \eqref{eqn:7} and initial condition \eqref{eqn:8} until time $t_n+\tau$. 

\textbf{(iv) Restriction:}
The average	
\begin{equation}\label{eqn:10}
U_j^{1,n}=\frac{1}{h}\int_{x_{j}-\frac{h}{2}}^{x_{j}+\frac{h}{2}}u_j(\xi,t_n+\tau)d\xi,\hspace{0.5cm} \forall j=1,...,N-1,
\end{equation}	
is computed which approximates the discretized macroscale variables, denoted by $U_j^{1,n}$ at position $x_j$ and at time $t_n+\tau$. 

\begin{figure}[h]
	\centering
	\advance\leftskip-3cm
	\advance\rightskip-3cm
	\begin{tikzpicture}
	\draw[ultra thick] (5,-2) -- (14,-2);
	
	\filldraw[thick, top color=white,bottom color=red!50!] (6,-1.98) rectangle node{} +(1,3);
	\filldraw[thick, top color=white,bottom color=red!50!] (9,-1.98) rectangle node{} +(1,3);
	\filldraw[thick, top color=white,bottom color=red!50!] (12,-1.98) rectangle node{} +(1,3);
	
	\filldraw[color=blue!60, fill=green!60, very thick](6.5,-1.98) circle (0.1);
	\filldraw[color=blue!60, fill=green!60, very thick](9.5,-1.98) circle (0.1);
	\filldraw[color=blue!60, fill=green!60, very thick](12.5,-1.98) circle (0.1);
	
	\filldraw[color=red!60, fill=yellow!60, very thick](6.5,1) circle (0.1);
	\filldraw[color=red!60, fill=yellow!60, very thick](9.5,1) circle (0.1);
	\filldraw[color=red!60, fill=yellow!60, very thick](12.5,1) circle (0.1);
	
	\draw[very thick] (6.5,-2.3)  node {$x_{j-1}$};
	\draw[very thick] (9.5,-2.3)  node {$x_j$};
	\draw[very thick] (12.5,-2.3)  node {$x_{j+1}$};
	
	\draw[very thick] (6.5,-0.1) node {\textbf{Tooth}};
	\draw[very thick] (9.5,-0.1) node {\textbf{Tooth}};
	\draw[very thick] (12.5,-0.1) node {\textbf{Tooth}};
	\draw[very thick] (8,-0.1) node {\textbf{Gap}};
	\draw[very thick] (11,-0.1) node {\textbf{Gap}};
	
	\draw[ultra thick, ->, blue] (9,-2) arc (95:150:25pt);
	\draw [ultra thick, ->, blue] (10,-2) arc (80:30:25pt);
	
	%\draw[ultra thick, ->] (1,0) arc (95:150:5pt);
	%\draw [ultra thick, ->] (1.1,0) arc (80:30:5pt);

	\draw[very thick] (8.5,-2.6)  node {$x_j-\frac{h}{2}$};
	\draw[very thick] (10.5,-2.6)  node {$x_j+\frac{h}{2}$};
	%\draw[very thick] (2,-2.6)  node {};
	%\draw[very thick] (2,+2)  node {};
	
	\draw[ultra thick, green] (4,-2)  node {Lift (IC): };
	%\draw[ultra thick] (0.1,-0.8)  node {Lifting(IC):};
	%\draw[very thick] (0.5,-0.15)  node {$u_x|_{X_j+\frac{h}{2}}=s_j^+$};
	\draw[ultra thick,green, ->] (5,-2) -- (5.8,-2);
	%\draw[ultra thick] (-0.15,0.45) -- (-0.15,0.503);
	
	\draw[ultra thick,blue, ->] (5,-1.6) -- (5,0.5);
	\draw[ultra thick, blue] (4,-0.5)  node {Evolution};
	
	\draw[ultra thick,red, ->] (6.45,1.05) arc (55:120:40pt);
	\draw[ultra thick,red] (4,1)  node {Restriction:};

	\end{tikzpicture}
	\caption{A schematic description of the gap-tooth timestepper (GTT).} \label{fig:1}
\end{figure}
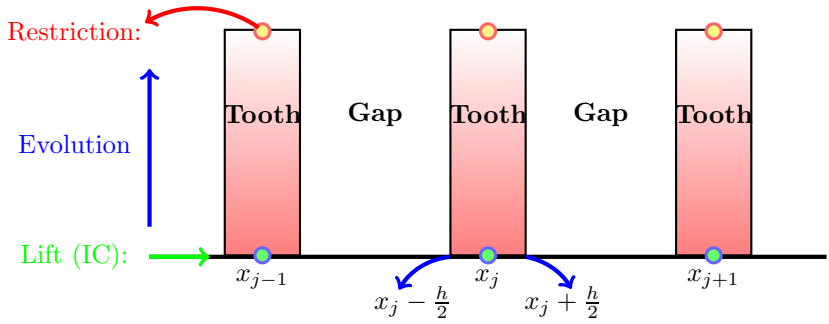

\textbf{Coarse Projective Integration:}

\textbf{(v) Short Time Steps:} The GTT is repeated for $k$ times, to reach at time $t_n+k\tau \hspace{0.2cm}(\ll t_n+\Delta t=t_{n+1})$ to find the value of $U_j^{k,n}$, \hspace{0.1cm} $\forall$ $j=1,...,N-1$. One more GTT is computed to find $U_j^{k+1,n}$ at time $t_n+(k+1)\tau$, which will help to evaluate an approximate value of the time derivative for the macroscopic variable, $\frac{\partial U}{\partial t}$ at $(x_j,t_n+k\tau)$.	

The time derivative of the macro state in each tooth at time $t_n+k\tau$ is estimated as 

\begin{equation}\label{eqn:11}
F(U^{k,n}, x_j, t_n+k\tau)=\frac{U_j^{k+1, n}-U_j^{k, n}}{\tau}, \hspace{0.3cm} \forall j=1,...,N-1.
\end{equation}

\textbf{(vi) Extrapolation:} 
The above estimate is used within time-stepping method, like forward Euler or Runge-Kutta, to march U forward in a long time to reach the macroscopic time $t_{n+1}=t_n+\Delta t$. For forward Euler, 

\begin{equation}\label{eqn:12}
U_j^{0,n+1}=U_j^{k,n}+(\Delta t-k\tau)F(U^{k,n}, x_j, t_n+k\tau), \hspace{0.3cm} \forall j=1,...,N-1.
\end{equation}

The whole procedure (from (i) to (vi)) is repeated in each long time steps $\Delta t$ until reach the final time. In this way, one can find the coarse solution using the microscopic problem, where the computation is done in only in a fraction of space and in a fraction of time.	

\begin{figure}[h!]
	\centering
	\begin{tikzpicture}[scale=8.5]
	\draw[ultra thick] (0,0) -- (1.2,0);
	\draw[ultra thick] (0,0.76) -- (1.2,0.76);
	
	\filldraw[thick, top color=white,bottom color=red!50!] (0.1,0) rectangle node{} +(0.1,0.05);
	\filldraw[thick, top color=white,bottom color=red!50!] (0.55,0) rectangle node{} +(0.1,0.05);
	\filldraw[thick, top color=white,bottom color=red!50!] (1,0) rectangle node{} +(0.1,0.05);
	
	\filldraw[thick, top color=white,bottom color=yellow!50!] (0.1,0.05) rectangle node{} +(0.1,0.05);
	\filldraw[thick, top color=white,bottom color=yellow!50!] (0.55,0.05) rectangle node{} +(0.1,0.05);
	\filldraw[thick, top color=white,bottom color=yellow!50!] (1,0.05) rectangle node{} +(0.1,0.05);
	
	\filldraw[thick, top color=white,bottom color=blue!50!] (0.1,0.1) rectangle node{} +(0.1,0.05);
	\filldraw[thick, top color=white,bottom color=blue!50!] (0.55,0.1) rectangle node{} +(0.1,0.05);
	\filldraw[thick, top color=white,bottom color=blue!50!] (1,0.1) rectangle node{} +(0.1,0.05);
	
	\filldraw[thick, top color=white,bottom color=magenta!50!] (0.1,0.15) rectangle node{} +(0.1,0.05);
	\filldraw[thick, top color=white,bottom color=magenta!50!] (0.55,0.15) rectangle node{} +(0.1,0.05);
	\filldraw[thick, top color=white,bottom color=magenta!50!] (1,0.15) rectangle node{} +(0.1,0.05);
	
	\filldraw[thick, top color=white,bottom color=violet!50!] (0.1,0.2) rectangle node{} +(0.1,0.05);
	\filldraw[thick, top color=white,bottom color=violet!50!] (0.55,0.2) rectangle node{} +(0.1,0.05);
	\filldraw[thick, top color=white,bottom color=violet!50!] (1,0.2) rectangle node{} +(0.1,0.05);
	
	\filldraw[thick, top color=white,bottom color=teal!50!] (0.1,0.25) rectangle node{} +(0.1,0.05);
	\filldraw[thick, top color=white,bottom color=teal!50!] (0.55,0.25) rectangle node{} +(0.1,0.05);
	\filldraw[thick, top color=white,bottom color=teal!50!] (1,0.25) rectangle node{} +(0.1,0.05);
	
	\filldraw[thick, top color=white,bottom color=brown!50!] (0.1,0.3) rectangle node{} +(0.1,0.05);
	\filldraw[thick, top color=white,bottom color=brown!50!] (0.55,0.3) rectangle node{} +(0.1,0.05);
	\filldraw[thick, top color=white,bottom color=brown!50!] (1,0.3) rectangle node{} +(0.1,0.05);
	
	\draw[ultra thick] (0.15,0.3) -- (0.15,0.75);
	\draw[ultra thick] (0.6,0.3) -- (0.6,0.75);
	\draw[ultra thick] (1.05,0.3) -- (1.05,0.75);
	
	\filldraw[color=blue!60, fill=green!60, very thick](0.15,0) circle (0.01);
	\filldraw[color=blue!60, fill=green!60, very thick](0.6,0) circle (0.01);
	\filldraw[color=blue!60, fill=green!60, very thick](1.05,0) circle (0.01);
	
	\draw[very thick] (0.15,-0.05)  node {$x_{j-1}$};
	\draw[very thick] (0.6,-0.05)  node {$x_j$};
	\draw[very thick] (1.05,-0.05)  node {$x_{j+1}$};
	
	\filldraw[color=red!60, fill=red!5, very thick](0.15,0.05) circle (0.01);
	\filldraw[color=red!60, fill=red!5, very thick](0.6,0.05) circle (0.01);
	\filldraw[color=red!60, fill=red!5, very thick](1.05,0.05) circle (0.01);
	
	\filldraw[color=red!60, fill=red!5, very thick](0.15,0.1) circle (0.01);
	\filldraw[color=red!60, fill=red!5, very thick](0.6,0.1) circle (0.01);
	\filldraw[color=red!60, fill=red!5, very thick](1.05,0.1) circle (0.01);
	
	\filldraw[color=red!60, fill=red!5, very thick](0.15,0.15) circle (0.01);
	\filldraw[color=red!60, fill=red!5, very thick](0.6,0.15) circle (0.01);
	\filldraw[color=red!60, fill=red!5, very thick](1.05,0.15) circle (0.01);
	
	\filldraw[color=red!60, fill=red!5, very thick](0.15,0.2) circle (0.01);
	\filldraw[color=red!60, fill=red!5, very thick](0.6,0.2) circle (0.01);
	\filldraw[color=red!60, fill=red!5, very thick](1.05,0.2) circle (0.01);
	
	\filldraw[color=red!60, fill=red!5, very thick](0.15,0.25) circle (0.01);
	\filldraw[color=red!60, fill=red!5, very thick](0.6,0.25) circle (0.01);
	\filldraw[color=red!60, fill=red!5, very thick](1.05,0.25) circle (0.01);
	
	\filldraw[color=blue!60, fill=green!60, very thick] (0.15,0.3) circle (0.01);
	\filldraw[color=blue!60, fill=green!60, very thick] (0.6,0.3) circle (0.01);
	\filldraw[color=blue!60, fill=green!60, very thick] (1.05,0.3) circle (0.01);
	
	\filldraw[color=red!60, fill=red!5, very thick] (0.15,0.35) circle (0.01);
	\filldraw[color=red!60, fill=red!5, very thick] (0.6,0.35) circle (0.01);
	\filldraw[color=red!60, fill=red!5, very thick] (1.05,0.35) circle (0.01);
	
	\filldraw[color=blue!60, fill=green!60, very thick](0.15,0.76) circle (0.01);
	\filldraw[color=blue!60, fill=green!60, very thick](0.6,0.76) circle (0.01);
	\filldraw[color=blue!60, fill=green!60, very thick](1.05,0.76) circle (0.01);
	
	\draw[very thick] (-0.07,0)  node {$t=t_n$};
	\draw[very thick] (-0.1,0.77)  node {$t=t_{n+1}$};
	\draw[very thick] (-0.06,0.73)  node {$=t_n+\Delta t$};
	
	\draw[ultra thick] (0.15,0.225)  node {Tooth};
	\draw[ultra thick] (0.38,0.22)  node {Gap};
	\draw[->, ultra thick] (0.33,0.22) -- (0.25,0.22);
	\draw[->, ultra thick] (0.42,0.22) -- (0.50,0.22);
	\draw[ultra thick] (0.6,0.225)  node {Tooth};
	\draw[ultra thick] (0.82,0.22)  node {Gap};
	\draw[->, ultra thick] (0.775,0.22) -- (0.69,0.22);
	\draw[->, ultra thick] (0.87,0.22) -- (0.95,0.22);
	\draw[ultra thick] (1.05,0.225)  node {Tooth};
	
	\draw [very thick, decorate,decoration={brace,amplitude=5pt},xshift=5pt,yshift=3pt]
	(-0.1,0.2) -- (-0.1,0.65) node [blue,midway,xshift=0cm] 
	{\footnotesize };
	\draw [very thick, decorate,decoration={brace,amplitude=5pt},xshift=5pt,yshift=3pt]
	(-0.13,-0.1) -- (-0.13,0.24) node [blue,midway,yshift=0cm] 
	{\footnotesize };
	
	\draw[ultra thick] (-0.075,0.53)  node {Extrapolation};
	\draw[ultra thick] (-0.03,0.17)  node {GTS};
	
	\draw[very thick] (0.25,0.3)  node {$U_{j-1}^{k,n}$};
	\draw[very thick] (0.7,0.3)  node {$U_j^{k,n}$};
	\draw[very thick] (1.15,0.3)  node {$U_{j+1}^{k,n}$};
	
	\draw[very thick] (0.15,0.8)  node {$U_{j-1}^{0,n+1}$};
	\draw[very thick] (0.6,0.8)  node {$U_j^{0,n+1}$};
	\draw[very thick] (1.05,0.8)  node {$U_{j+1}^{0,n+1}$};
	
	\draw[very thick] (0.25,0.08)  node {$U_{j-1}^{m,n}$};
	\draw[very thick] (0.7,0.08)  node {$U_j^{m,n}$};
	\draw[very thick] (1.15,0.08)  node {$U_{j+1}^{m,n}$};
	
	\draw[very thick] (0.27,0.15)  node {$U_{j-1}^{m+1,n}$};
	\draw[very thick] (0.72,0.15)  node {$U_j^{m+1,n}$};
	\draw[very thick] (1.17,0.15)  node {$U_{j+1}^{m+1,n}$};
	
	\end{tikzpicture}
	\caption{A schematic description of the usual patch dynamics scheme.} \label{fig:2}
\end{figure}
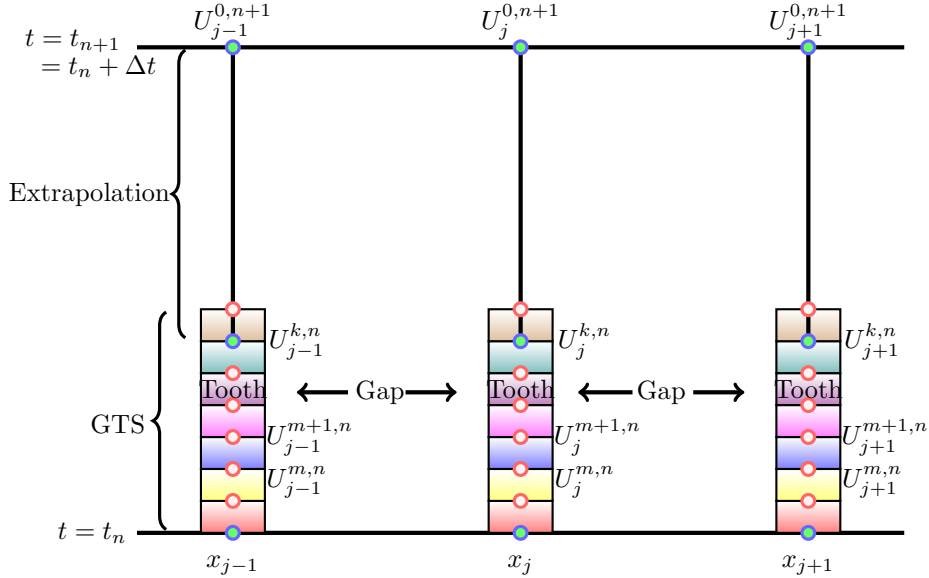

In figure \ref{fig:2} , $X_{j-1}, X_j$ and $X_{j+1}$ are the coarse grid points. Suppose $k$ number of GTTs (from bottom red box to teal box) are implemented to find $U_j^{k,n}$, subsequently one extra GTT (brown box at the top) is implemented to estimate the time-derivative $\frac{\partial U}{\partial t}$ at $(x_j,t_n+k\tau)$, and finally, extrapolation (black vertical line) is applied here. Suppose $\tau$ is the evolution time for each GTT. Then total microscopic computation is done for time $(k+1)\tau$, and extrapolation is done for time ($\Delta t$-$k\tau$) to perform a long time integration on time $\Delta t$.

\section{Problem Discussion}
After a series of experiments, it has been observed that the extrapolation time step for long time integration in UPD scheme is an important factor for the accuracy and convergence of the solution. Due to this factor, the solution always contains a significant error. So, it is very important to choose a correct extrapolation step size to reduce the error in the solution to achieve a better accuracy. It is obvious that if the extrapolation time step size is reduced, a better accuracy in the solution can be achieved. This idea can be implimented in two ways. Firstly, reduce the time-step $\Delta t$ by maintaining the same `$k$' times GTTs such that $k\tau \ll \Delta t$ for time integration $\Delta t$.  That means when $\Delta t$ is reduced, effectively the size of extrapolation step $\Delta t-k\tau (> 0)$ is also reduced. In another way, total `$k$' number of GTTs can be distributed uniformly or non-uniformly within the long time step $\Delta t$. This kind of distribution has some advantages over the UPD scheme. These distributions provide the flexibility to incorporate three different time scales, namely micro, meso and macro time scales, to accurately capture the physical behaviour over time. The macroscale solver utilizes a mesoscale time step finer than the macro time scale to determind the system level behaviour. This new multiscale integrator can be applied to handle stiff and highly oscillatory problems more efficiently. The distribution of GTTs can be done in two ways, namely generalized patch dynamics (GPD) schemes of type-I and type-II, which are discussed in the next section. GPD scheme of type-I is a `GTT \textbf{-} extrapolation' based modelling (fig- \ref{fig:3}) and GPD scheme of type-II is a `GTT \textbf{-} extrapolation \textbf{-} GTT' based modelling (fig- \ref{fig:4}).

\subsection{Generalized Patch Dynamics (GPD) Scheme of Type-I}

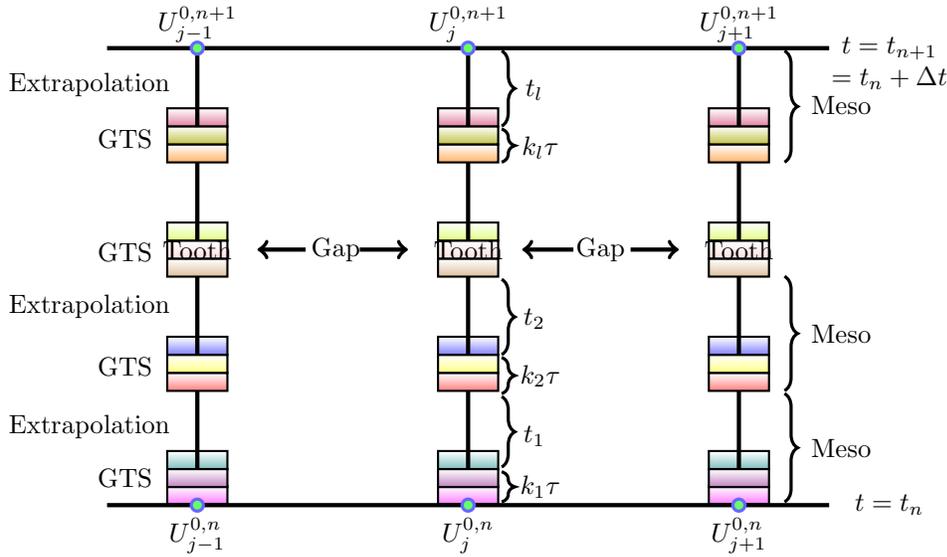
\begin{figure}[h!]
	\begin{tikzpicture}[scale=8]
	\draw[ultra thick] (0,0) -- (1.2,0);
	\draw[ultra thick] (0,0.76) -- (1.2,0.76);
	
	\filldraw[thick, top color=white,bottom color=magenta!50!] (0.1,0.002) rectangle node{} +(0.1,0.03);
	\filldraw[thick, top color=white,bottom color=magenta!50!] (0.55,0.002) rectangle node{} +(0.1,0.03);
	\filldraw[thick, top color=white,bottom color=magenta!50!] (1,0.002) rectangle node{} +(0.1,0.03);
	
	\filldraw[thick, top color=white,bottom color=violet!50!] (0.1,0.03) rectangle node{} +(0.1,0.03);
	\filldraw[thick, top color=white,bottom color=violet!50!] (0.55,0.03) rectangle node{} +(0.1,0.03);
	\filldraw[thick, top color=white,bottom color=violet!50!] (1,0.03) rectangle node{} +(0.1,0.03);
	
	\filldraw[thick, top color=white,bottom color=teal!50!] (0.1,0.06) rectangle node{} +(0.1,0.03);
	\filldraw[thick, top color=white,bottom color=teal!50!] (0.55,0.06) rectangle node{} +(0.1,0.03);
	\filldraw[thick, top color=white,bottom color=teal!50!] (1,0.06) rectangle node{} +(0.1,0.03);
	
	\draw[ultra thick] (0.15,0.06) -- (0.15,0.19);
	\draw[ultra thick] (0.6,0.06) -- (0.6,0.19);
	\draw[ultra thick] (1.05,0.06) -- (1.05,0.19);
	
	\filldraw[thick, top color=white,bottom color=red!50!] (0.1,0.19) rectangle node{} +(0.1,0.03);
	\filldraw[thick, top color=white,bottom color=red!50!] (0.55,0.19) rectangle node{} +(0.1,0.03);
	\filldraw[thick, top color=white,bottom color=red!50!] (1,0.19) rectangle node{} +(0.1,0.03);
	
	\filldraw[thick, top color=white,bottom color=yellow!50!] (0.1,0.22) rectangle node{} +(0.1,0.03);
	\filldraw[thick, top color=white,bottom color=yellow!50!] (0.55,0.22) rectangle node{} +(0.1,0.03);
	\filldraw[thick, top color=white,bottom color=yellow!50!] (1,0.22) rectangle node{} +(0.1,0.03);
	
	\filldraw[thick, top color=white,bottom color=blue!50!] (0.1,0.25) rectangle node{} +(0.1,0.03);
	\filldraw[thick, top color=white,bottom color=blue!50!] (0.55,0.25) rectangle node{} +(0.1,0.03);
	\filldraw[thick, top color=white,bottom color=blue!50!] (1,0.25) rectangle node{} +(0.1,0.03);
	
	\draw[ultra thick] (0.15,0.25) -- (0.15,0.38);
	\draw[ultra thick] (0.6,0.25) -- (0.6,0.38);
	\draw[ultra thick] (1.05,0.25) -- (1.05,0.38);
	
	\filldraw[thick, top color=white,bottom color=brown!50!] (0.1,0.38) rectangle node{} +(0.1,0.03);
	\filldraw[thick, top color=white,bottom color=brown!50!] (0.55,0.38) rectangle node{} +(0.1,0.03);
	\filldraw[thick, top color=white,bottom color=brown!50!] (1,0.38) rectangle node{} +(0.1,0.03);
	
	\filldraw[thick, top color=white,bottom color=pink!50!] (0.1,0.41) rectangle node{} +(0.1,0.03);
	\filldraw[thick, top color=white,bottom color=pink!50!] (0.55,0.41) rectangle node{} +(0.1,0.03);
	\filldraw[thick, top color=white,bottom color=pink!50!] (1,0.41) rectangle node{} +(0.1,0.03);
	
	\filldraw[thick, top color=white,bottom color=lime!50!] (0.1,0.44) rectangle node{} +(0.1,0.03);
	\filldraw[thick, top color=white,bottom color=lime!50!] (0.55,0.44) rectangle node{} +(0.1,0.03);
	\filldraw[thick, top color=white,bottom color=lime!50!] (1,0.44) rectangle node{} +(0.1,0.03);
	
	\draw[ultra thick] (0.15,0.44) -- (0.15,0.57);
	\draw[ultra thick] (0.6,0.44) -- (0.6,0.57);
	\draw[ultra thick] (1.05,0.44) -- (1.05,0.57);
	
	\filldraw[thick, top color=white,bottom color=orange!60!] (0.1,0.57) rectangle node{} +(0.1,0.03);
	\filldraw[thick, top color=white,bottom color=orange!60!] (0.55,0.57) rectangle node{} +(0.1,0.03);
	\filldraw[thick, top color=white,bottom color=orange!60!] (1,0.57) rectangle node{} +(0.1,0.03);
	
	\filldraw[thick, top color=white,bottom color=olive!50!] (0.1,0.6) rectangle node{} +(0.1,0.03);
	\filldraw[thick, top color=white,bottom color=olive!50!] (0.55,0.6) rectangle node{} +(0.1,0.03);
	\filldraw[thick, top color=white,bottom color=olive!50!] (1,0.6) rectangle node{} +(0.1,0.03);
	
	\filldraw[thick, top color=white,bottom color=purple!50!] (0.1,0.63) rectangle node{} +(0.1,0.03);
	\filldraw[thick, top color=white,bottom color=purple!50!] (0.55,0.63) rectangle node{} +(0.1,0.03);
	\filldraw[thick, top color=white,bottom color=purple!50!] (1,0.63) rectangle node{} +(0.1,0.03);
	
	\draw[ultra thick] (0.15,0.63) -- (0.15,0.759);
	\draw[ultra thick] (0.6,0.63) -- (0.6,0.759);
	\draw[ultra thick] (1.05,0.63) -- (1.05,0.759);
	
	\filldraw[color=blue!60, fill=green!60, very thick](0.15,0) circle (0.01);
	\filldraw[color=blue!60, fill=green!60, very thick](0.6,0) circle (0.01);
	\filldraw[color=blue!60, fill=green!60, very thick](1.05,0) circle (0.01);
	
	\filldraw[color=blue!60, fill=green!60, very thick](0.15,0.76) circle (0.01);
	\filldraw[color=blue!60, fill=green!60, very thick](0.6,0.76) circle (0.01);
	\filldraw[color=blue!60, fill=green!60, very thick](1.05,0.76) circle (0.01);
	
	\draw[very thick] (1.3,0)  node {$t=t_n$};
	\draw[very thick] (1.3,0.76)  node {$t=t_{n+1}$};
	\draw[very thick] (1.3,0.71)  node {$=t_n+\Delta t$};
	
	\draw[ultra thick] (0.15,0.425)  node {Tooth};
	\draw[ultra thick] (0.38,0.425)  node {Gap};
	\draw[->, ultra thick] (0.33,0.425) -- (0.25,0.425);
	\draw[->, ultra thick] (0.42,0.425) -- (0.50,0.425);
	\draw[ultra thick] (0.6,0.425)  node {Tooth};
	\draw[ultra thick] (0.82,0.425)  node {Gap};
	\draw[->, ultra thick] (0.775,0.425) -- (0.69,0.425);
	\draw[->, ultra thick] (0.87,0.425) -- (0.95,0.425);
	\draw[ultra thick] (1.05,0.425)  node {Tooth};
	
	\draw[ultra thick] (0.03,0.045)  node {GTS};
	\draw [very thick, decorate,decoration={brace,mirror, amplitude=5pt},xshift=5pt,yshift=3pt]
	(0.48,-0.1) -- (0.48,-0.05) node [blue,midway,yshift=0cm] 
	{\footnotesize };
	\draw[ultra thick] (0.72,0.03)  node {$k_1\tau$};
	
	\draw[ultra thick] (-0.03,0.13)  node {Extrapolation};
	\draw [very thick, decorate,decoration={brace,mirror,amplitude=5pt},xshift=5pt,yshift=3pt]
	(0.48,-0.042) -- (0.48,0.075) node [blue,midway,yshift=0cm] 
	{\footnotesize };
	\draw[ultra thick] (0.71,0.115)  node {$t_1$};
	\draw [very thick, decorate,decoration={brace,mirror,amplitude=5pt},xshift=5pt,yshift=3pt]
	(0.95,-0.1) -- (0.95,0.08) node [blue,midway,yshift=0cm] 
	{\footnotesize };
	\draw[ultra thick] (1.22,0.095)  node {Meso};
	%----------------------------
	\draw[ultra thick] (0.03,0.23)  node {GTS};
	\draw [very thick, decorate,decoration={brace,mirror,amplitude=5pt},xshift=5pt,yshift=3pt]
	(0.48,0.08) -- (0.48,0.14) node [blue,midway,yshift=0cm] 
	{\footnotesize };
	\draw[ultra thick] (0.72,0.215)  node {$k_2\tau$};
	
	\draw[ultra thick] (-0.03,0.33)  node {Extrapolation};
	\draw [very thick, decorate,decoration={brace,mirror,amplitude=5pt},xshift=5pt,yshift=3pt]
	(0.48,0.145) -- (0.48,0.27) node [blue,midway,yshift=0cm] 
	{\footnotesize };
	\draw[ultra thick] (0.71,0.31)  node {$t_2$};
	\draw [very thick, decorate,decoration={brace,mirror,amplitude=5pt},xshift=5pt,yshift=3pt]
	(0.95,0.085) -- (0.95,0.275) node [blue,midway,yshift=0cm] 
	{\footnotesize };
	\draw[ultra thick] (1.22,0.285)  node {Meso};
	%---------------------------
	\draw[ultra thick] (0.03,0.42)  node {GTS};
	%\draw[ultra thick] (-0.03,0.51)  node {Extrapolation};
	
	\draw[ultra thick] (0.03,0.61)  node {GTS};
	\draw [very thick, decorate,decoration={brace,mirror,amplitude=5pt},xshift=5pt,yshift=3pt]
	(0.48,0.465) -- (0.48,0.52) node [blue,midway,yshift=0cm] 
	{\footnotesize };
	\draw[ultra thick] (0.72,0.595)  node {$k_l\tau$};
	
	\draw[ultra thick] (-0.03,0.7)  node {Extrapolation};
	\draw [very thick, decorate,decoration={brace,mirror,amplitude=5pt},xshift=5pt,yshift=3pt]
	(0.48,0.525) -- (0.48,0.65) node [blue,midway,yshift=0cm] 
	{\footnotesize };
	\draw[ultra thick] (0.71,0.69)  node {$t_l$};
	\draw [very thick, decorate,decoration={brace,mirror,amplitude=5pt},xshift=5pt,yshift=3pt]
	(0.95,0.465) -- (0.95,0.65) node [blue,midway,yshift=0cm] 
	{\footnotesize };
	\draw[ultra thick] (1.22,0.665)  node {Meso};
	%-----------------------
	\draw[very thick] (0.15,-0.05)  node {$U_{j-1}^{0,n}$};
	\draw[very thick] (0.6,-0.05)  node {$U_j^{0,n}$};
	\draw[very thick] (1.05,-0.05)  node {$U_{j+1}^{0,n}$};
	
	\draw[very thick] (0.15,0.8)  node {$U_{j-1}^{0,n+1}$};
	\draw[very thick] (0.6,0.8)  node {$U_j^{0,n+1}$};
	\draw[very thick] (1.05,0.8)  node {$U_{j+1}^{0,n+1}$};
	
	%	\draw[ultra thick] (1.15,0) -- (1.15,0.76);
	%	\draw[ultra thick] (0.6,0.44) -- (0.6,0.57);
	%	\draw[ultra thick] (1.05,0.44) -- (1.05,0.57);
	
	\end{tikzpicture}
	\caption{A schematic description of the GPD scheme of Type-I. } \label{fig:3}
\end{figure}

In order to progress for a long macroscopic time $\Delta t$ in UPD scheme, first `$k$' number of GTTs are applied, subsequently the solution is extrapolated for the rest of the time by \eqref{eqn:12}. Suppose that `$k$' number of GTTs are distributed into `$l$' different parts separated from each other in time $\Delta t$ and the corresponding number of GTTs in `$l$' different parts are $k_1, k_2,..., k_l$ such that 

\begin{equation}\label{eqn:13}
k_1+k_2+k_3+ ... +k_l=k,
\end{equation}
in time $\Delta t$ as shown in figure \ref{fig:3}. Here, the GTTs are evaluated in a similar way as discussed in section 3 (from (i) to (iv)). Similarly, at the end of each of the group of GTTs, it is needed to compute one more GTT to find the temporal derivative term as \eqref{eqn:11} for long time extrapolation in meso time step. Let, the corresponding extrapolation time be $t_i$ after $k_i$ number of GTTs such that $t_i<\Delta t$, for all $i=1, 2,..., l$. Hence, $\l$ number of extrapolations are applied for time $t_1, t_2, ... ,t_l$ in each macro time step $\Delta t$ such that,

\begin{equation}\label{eqn:14}
t_1+t_2+t_3+ ... +t_l=\Delta t-(k_1+k_2+k_3+ ... +k_l)\tau=\Delta t-k\tau.
\end{equation}

A schematic description of the GPD scheme of type-I is shown in figure \ref{fig:3} which is a `GTT-extrapolation' based modelling.

One has the freedom to choose uniform or non-uniform distribution of GTTs as well as extrapolation times based on the physical problem. For uniform distribution of GTTs, all $k_i$'s are equal, that is $k_1=k_2=...=k_l=\frac{k}{l}$, (from \eqref{eqn:13}). Similarly, for uniform extrapolation times, all $t_i$'s are also equal, that is $t_1=t_2=...=t_l=\frac{\Delta t-k\tau}{l}$, (from \eqref{eqn:14}).
For non-uniform GPD scheme of type-I, one can employ a variable number of microsimulators and variable mesoscopic time steps.

\subsection{Generalized Patch Dynamics (GPD) Scheme of Type-II}
GPD scheme of type-II is a `GTT \textbf{-} extrapolation \textbf{-} GTT' based modelling. The GPD scheme of type-II is a modification of the GPD scheme of type-I. Here, microsimulations and extrapolations are alternatively implemented, following the structure of the GPD scheme of type-I. However, the final extrapolation is executed just before the time $t_{n+1}$, and the last group of the GTTs is implemented at the end of the macro time step $\Delta t$. A similar idea can be observed in modified projective integration, named PI2, proposed by Maclean et al. \cite{2015_maclean_convergence}. Suppose the total `$k$' number of GTTs are distributed into `$l$' different parts in each time step $\Delta t$ and the corresponding number of GTTs in `$l$' different parts are $k_1, k_2,..., k_l$ such that 

\begin{figure}[h!]	
	\begin{tikzpicture}[scale=8]
	\draw[ultra thick] (0,0) -- (1.2,0);
	\draw[ultra thick] (0,0.661) -- (1.2,0.661);
	
	\filldraw[thick, top color=white,bottom color=magenta!50!] (0.1,0.002) rectangle node{} +(0.1,0.03);
	\filldraw[thick, top color=white,bottom color=magenta!50!] (0.55,0.002) rectangle node{} +(0.1,0.03);
	\filldraw[thick, top color=white,bottom color=magenta!50!] (1,0.002) rectangle node{} +(0.1,0.03);
	
	\filldraw[thick, top color=white,bottom color=violet!50!] (0.1,0.03) rectangle node{} +(0.1,0.03);
	\filldraw[thick, top color=white,bottom color=violet!50!] (0.55,0.03) rectangle node{} +(0.1,0.03);
	\filldraw[thick, top color=white,bottom color=violet!50!] (1,0.03) rectangle node{} +(0.1,0.03);
	
	\filldraw[thick, top color=white,bottom color=teal!50!] (0.1,0.06) rectangle node{} +(0.1,0.03);
	\filldraw[thick, top color=white,bottom color=teal!50!] (0.55,0.06) rectangle node{} +(0.1,0.03);
	\filldraw[thick, top color=white,bottom color=teal!50!] (1,0.06) rectangle node{} +(0.1,0.03);
	
	\draw[ultra thick] (0.15,0.06) -- (0.15,0.19);
	\draw[ultra thick] (0.6,0.06) -- (0.6,0.19);
	\draw[ultra thick] (1.05,0.06) -- (1.05,0.19);
	
	\filldraw[thick, top color=white,bottom color=red!50!] (0.1,0.19) rectangle node{} +(0.1,0.03);
	\filldraw[thick, top color=white,bottom color=red!50!] (0.55,0.19) rectangle node{} +(0.1,0.03);
	\filldraw[thick, top color=white,bottom color=red!50!] (1,0.19) rectangle node{} +(0.1,0.03);
	
	\filldraw[thick, top color=white,bottom color=yellow!50!] (0.1,0.22) rectangle node{} +(0.1,0.03);
	\filldraw[thick, top color=white,bottom color=yellow!50!] (0.55,0.22) rectangle node{} +(0.1,0.03);
	\filldraw[thick, top color=white,bottom color=yellow!50!] (1,0.22) rectangle node{} +(0.1,0.03);
	
	\filldraw[thick, top color=white,bottom color=blue!50!] (0.1,0.25) rectangle node{} +(0.1,0.03);
	\filldraw[thick, top color=white,bottom color=blue!50!] (0.55,0.25) rectangle node{} +(0.1,0.03);
	\filldraw[thick, top color=white,bottom color=blue!50!] (1,0.25) rectangle node{} +(0.1,0.03);
	
	\draw[ultra thick] (0.15,0.25) -- (0.15,0.38);
	\draw[ultra thick] (0.6,0.25) -- (0.6,0.38);
	\draw[ultra thick] (1.05,0.25) -- (1.05,0.38);
	
	\filldraw[thick, top color=white,bottom color=brown!50!] (0.1,0.38) rectangle node{} +(0.1,0.03);
	\filldraw[thick, top color=white,bottom color=brown!50!] (0.55,0.38) rectangle node{} +(0.1,0.03);
	\filldraw[thick, top color=white,bottom color=brown!50!] (1,0.38) rectangle node{} +(0.1,0.03);
	
	\filldraw[thick, top color=white,bottom color=pink!50!] (0.1,0.41) rectangle node{} +(0.1,0.03);
	\filldraw[thick, top color=white,bottom color=pink!50!] (0.55,0.41) rectangle node{} +(0.1,0.03);
	\filldraw[thick, top color=white,bottom color=pink!50!] (1,0.41) rectangle node{} +(0.1,0.03);
	
	\filldraw[thick, top color=white,bottom color=lime!50!] (0.1,0.44) rectangle node{} +(0.1,0.03);
	\filldraw[thick, top color=white,bottom color=lime!50!] (0.55,0.44) rectangle node{} +(0.1,0.03);
	\filldraw[thick, top color=white,bottom color=lime!50!] (1,0.44) rectangle node{} +(0.1,0.03);
	
	\draw[ultra thick] (0.15,0.44) -- (0.15,0.57);
	\draw[ultra thick] (0.6,0.44) -- (0.6,0.57);
	\draw[ultra thick] (1.05,0.44) -- (1.05,0.57);
	
	\filldraw[thick, top color=white,bottom color=orange!60!] (0.1,0.57) rectangle node{} +(0.1,0.03);
	\filldraw[thick, top color=white,bottom color=orange!60!] (0.55,0.57) rectangle node{} +(0.1,0.03);
	\filldraw[thick, top color=white,bottom color=orange!60!] (1,0.57) rectangle node{} +(0.1,0.03);
	
	\filldraw[thick, top color=white,bottom color=olive!50!] (0.1,0.6) rectangle node{} +(0.1,0.03);
	\filldraw[thick, top color=white,bottom color=olive!50!] (0.55,0.6) rectangle node{} +(0.1,0.03);
	\filldraw[thick, top color=white,bottom color=olive!50!] (1,0.6) rectangle node{} +(0.1,0.03);
	
	\filldraw[thick, top color=white,bottom color=purple!50!] (0.1,0.63) rectangle node{} +(0.1,0.03);
	\filldraw[thick, top color=white,bottom color=purple!50!] (0.55,0.63) rectangle node{} +(0.1,0.03);
	\filldraw[thick, top color=white,bottom color=purple!50!] (1,0.63) rectangle node{} +(0.1,0.03);
	
	\filldraw[color=blue!60, fill=green!60, very thick](0.15,0) circle (0.01);
	\filldraw[color=blue!60, fill=green!60, very thick](0.6,0) circle (0.01);
	\filldraw[color=blue!60, fill=green!60, very thick](1.05,0) circle (0.01);
	
	\filldraw[color=blue!60, fill=green!60, very thick](0.15,0.66) circle (0.01);
	\filldraw[color=blue!60, fill=green!60, very thick](0.6,0.66) circle (0.01);
	\filldraw[color=blue!60, fill=green!60, very thick](1.05,0.66) circle (0.01);
	
	\draw[very thick] (1.28,0)  node {$t=t_n$};
	\draw[very thick] (1.3,0.66)  node {$t=t_{n+1}$};
	\draw[very thick] (1.3,0.61)  node {$=t_n+\Delta t$};
	
	\draw[ultra thick] (0.15,0.425)  node {Tooth};
	\draw[ultra thick] (0.38,0.425)  node {Gap};
	\draw[->, ultra thick] (0.33,0.425) -- (0.25,0.425);
	\draw[->, ultra thick] (0.42,0.425) -- (0.50,0.425);
	\draw[ultra thick] (0.6,0.425)  node {Tooth};
	\draw[ultra thick] (0.82,0.425)  node {Gap};
	\draw[->, ultra thick] (0.775,0.425) -- (0.69,0.425);
	\draw[->, ultra thick] (0.87,0.425) -- (0.95,0.425);
	\draw[ultra thick] (1.05,0.425)  node {Tooth};
	
	\draw[ultra thick] (0.03,0.045)  node {GTS};
	\draw [very thick, decorate,decoration={brace,mirror, amplitude=5pt},xshift=5pt,yshift=3pt]
	(0.48,-0.1) -- (0.48,-0.05) node [blue,midway,yshift=0cm] 
	{\footnotesize };
	\draw[ultra thick] (0.72,0.03)  node {$k_1\tau$};
	
	\draw[ultra thick] (-0.03,0.13)  node {Extrapolation};
	\draw [very thick, decorate,decoration={brace,mirror,amplitude=5pt},xshift=5pt,yshift=3pt]
	(0.48,-0.043) -- (0.48,0.075) node [blue,midway,yshift=0cm] 
	{\footnotesize };
	\draw[ultra thick] (0.71,0.12)  node {$t_1$};
	
	%----------------------
	\draw[ultra thick] (0.03,0.23)  node {GTS};
	\draw [very thick, decorate,decoration={brace,mirror,amplitude=5pt},xshift=5pt,yshift=3pt]
	(0.48,0.08) -- (0.48,0.14) node [blue,midway,yshift=0cm] 
	{\footnotesize };
	\draw[ultra thick] (0.72,0.215)  node {$k_2\tau$};
	
	\draw[ultra thick] (-0.03,0.33)  node {Extrapolation};
	\draw [very thick, decorate,decoration={brace,mirror,amplitude=5pt},xshift=5pt,yshift=3pt]
	(0.48,0.146) -- (0.48,0.27) node [blue,midway,yshift=0cm] 
	{\footnotesize };
	\draw[ultra thick] (0.71,0.31)  node {$t_2$};
	
	%-----------------------
	\draw[ultra thick] (0.03,0.42)  node {GTS};
	
	%\draw[ultra thick] (-0.03,0.51)  node {Extrapolation};
	\draw [very thick, decorate,decoration={brace,mirror,amplitude=5pt},xshift=5pt,yshift=3pt]
	(0.48,0.338) -- (0.48,0.46) node [blue,midway,yshift=0cm] 
	{\footnotesize };
	\draw[ultra thick] (0.73,0.5)  node {$t_{l-1}$};
	
	\draw[ultra thick] (0.03,0.61)  node {GTS};
	\draw [very thick, decorate,decoration={brace,mirror,amplitude=5pt},xshift=5pt,yshift=3pt]
	(0.48,0.465) -- (0.48,0.55) node [blue,midway,yshift=0cm] 
	{\footnotesize };
	\draw[ultra thick] (0.72,0.61)  node {$k_l\tau$};
	
	\draw[very thick] (0.15,-0.05)  node {$U_{j-1}^{0,n}$};
	\draw[very thick] (0.6,-0.05)  node {$U_j^{0,n}$};
	\draw[very thick] (1.05,-0.05)  node {$U_{j+1}^{0,n}$};
	
	\draw[very thick] (0.15,0.7)  node {$U_{j-1}^{0,n+1}$};
	\draw[very thick] (0.6,0.7)  node {$U_j^{0,n+1}$};
	\draw[very thick] (1.05,0.7)  node {$U_{j+1}^{0,n+1}$};
	
	\end{tikzpicture}
	\caption{A schematic description of the GPD scheme of Type-II.} \label{fig:4}
\end{figure}
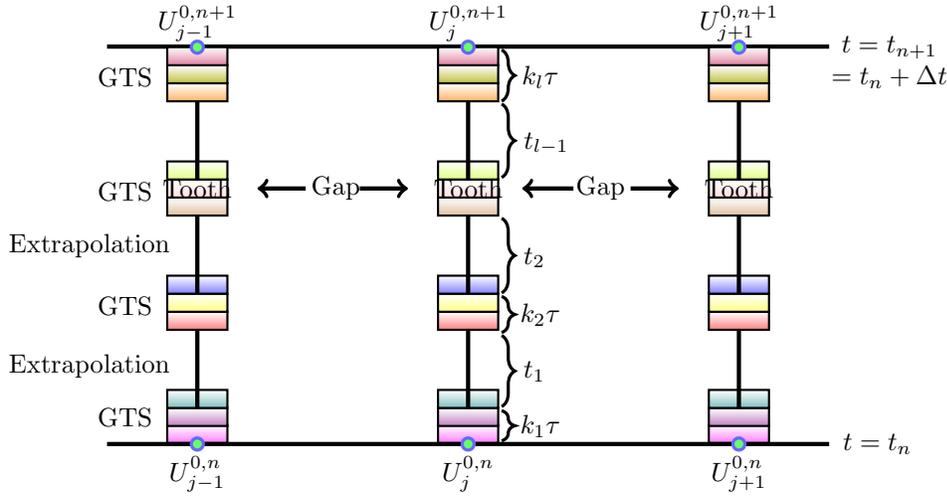

\begin{equation}\label{eqn:15}
k_1+k_2+k_3+ ... +k_l=k,
\end{equation}
in each time step $\Delta t$ as shown in figure \ref{fig:4}. The macro time step begins with $k_1$ GTTs and concludes with $k_l$ GTTs.  Here GTTs are evaluated in the same way as discussed in the section 3 (from (i) to (iv)). So at the end of each of the group of GTTs (except the last group), it is needed to compute one more GTT to find the temporal derivative term from \eqref{eqn:11} for long time extrapolation. As there is no requirement of extrapolating the solution after the last group of GTTs in macro time step $\Delta t$, no additional GTT is computed. Let, the corresponding extrapolation time be $t_i$ ($<\Delta t$) after completion of $k_i$ number of GTTs for all $i=1, 2,..., (l-1)$. Hence, $(\l-1)$ number of extrapolations are applied for time $t_1, t_2, ... ,t_{l-1}$ in each time step $\Delta t$ such that,

\begin{equation}\label{eqn:16}
t_1+t_2+t_3+ ... +t_{l-1}=\Delta t-(k_1+k_2+k_3+ ... +k_l)\tau=\Delta t-k\tau.
\end{equation}

In GPD scheme of type-II, also one can choose uniform or non-uniform distribution of GTTs as well as extrapolation times based on the physical problem. For uniform distribution of GTTs, all $k_i$'s are equal, that is $k_1=k_2=...=k_l=\frac{k}{l},$ (from \eqref{eqn:15}). Similarly, for uniform extrapolation times, all $t_i$'s are also equal, that is $t_1=t_2=...=t_{l-1}=\frac{\Delta t-k\tau}{l-1},$ (from \eqref{eqn:16}). However, the GPD scheme of type-II is a geneal concept that provides the flexibility to use a variable number of microsimulators and variable mesoscopic time steps.

%The GPD scheme is a generalized version of the UPD scheme.
%
%$\bullet$ The UPD scheme as proposed by Kevrekidis et al. \cite{2003_kevrekidis_equation} can be obtained from the GPD scheme of type-I by putting $l= 1$.
%
%$\bullet$ Similarly, the UPD scheme as proposed by Samaye et al. \cite{2006_samaey_patch,2020_arbabi_linking} can be obtained from the GPD scheme of type-I by putting $l= 1$ \& $k=0$.

The GPD schemes are expected to achieve a better accuracy due to the following reasons:

\begin{enumerate}
	\item \textbf{Reduced Extrapolation Time:}
	
	In comparison to the UPD scheme, the GPD schemes need a smaller extrapolation time in meso time steps compared to the macro time step. This reduction in time allows for more precise calculations in the GPD schemes.
	\item \textbf{Enhanced Number of GTTs:}
	
	In the GPD schemes, a few additional GTTs are simulated to compute the time derivatives (i.e., equation \eqref{eqn:11}) in the macro time steps for the outer integration. These extra microsimulations contribute to an improvement in the results.
	\item \textbf{Enhanced Spatial Bridging:}
	
	The GPD schemes should facilitate improved accuracy through enhanced communication between patches. Patches interact with each other through the interpolated macroscale and boundary conditions. Unlike the UPD scheme, where spatial gaps are bridged only ($k+1$) number of times at the bottom of the macro time interval [$t_n,t_n+\Delta t$), the GPD scheme of type-I and type-II, implements bridging the spatial gaps ($k+l$) and ($k+l-1$) number of times, respectively, at different time locations within the macro time interval [$t_n,t_n+\Delta t$). The increased frequency of bridging spatial gaps in the GPD schemes results in more effective tranformation and transfer of information compared to that in the UPD scheme.
\end{enumerate}

A fundamental underlying concept in all three approaches (systematic up-scaling, HMM and equation-free approach) is that the local relaxation time for the microscopic process is much shorter than the macroscopic time scale of the system \cite{2009_weinan_general}. In other words, the microscopic processes are in equilibrium with the macroscopic state of the system in the UPD scheme. However, in the GPD scheme, the microscopic processes are in equilibrium with the mesoscopic state of the system. `$k_i$' must be large enough to allow the micro-state to attain its quasi-equilibrium state. If $\tau_\epsilon$ be the relaxation time for the microscopic model, then $\tau_\epsilon<<k_i\tau$.

The most effective approach to achieving better accuracy with less computational time compared to the UPD scheme, one should adopt to taking fewer macro time steps with a finer distribution of GTTs (i.e., more meso time steps) while still adhering to the relaxation time limit. 

From the observations, one can easily realize that the
GPD scheme is a generalized version of UPD scheme:

$\bullet$ By putting $l=1$ in the uniform GPD scheme of type-I, one can obtain the UPD scheme as proposed by Kevrekidis et al. \cite{2003_kevrekidis_equation}.

$\bullet$ By putting $l=1$ \& $k=0$ in the uniform GPD scheme of type-I, one can obtain the UPD scheme as proposed by Samaye et al. \cite{2006_samaey_patch,2009_kevrekidis_equation,2009_samaey_equation,2020_arbabi_linking}.

\subsection{Projective Integration Version of the GPD Scheme}
If the model is independent of the spatial variables, the evolution equation \eqref{eqn:1} transforms to a system of ODEs dependent on temporal variables. In order to solve such problems, projective integration (PI) could be a better option. This article introduces the projective integration version of the GPD scheme (PI-GPD), a generalized version of the PI  and PI version of UPD (PI-UPD). Table \ref{table:model_PI} presents these approaches.

\begin{table}[h!]
	\begin{center}
		\begin{tabular}{|p{3cm}|p{3.5cm}|p{5cm}|}
			\hline
			%			%\cline{2-3}
			%			&\multicolumn{2}{|c|}{Distribution of 24 number of gap-}\\
			%			&\multicolumn{2}{|c|}{tooth timesteppers in GPD scheme of}\\
			%\hline
			\textbf{Proposed scheme}&\textbf{Values of $l$ and / or $k$}&\textbf{The proposed scheme will be equivalent to}\\
			\hline
			\multirow{2}{*}{Uniform PI-GPD-I}
			&$l=1$&PI \cite{2003_gear_projective}\\\cline{2-3}
			&$l=1$ \& $k=0$ (i.e., $k_1=0$)&PI version of UPD \cite{2006_samaey_patch,2009_kevrekidis_equation,2020_arbabi_linking,2009_samaey_equation}\\
			\hline
		\end{tabular}
	\end{center}
	\caption{The PI-GPD scheme is a generalized version of PI and PI-UPD schemes.}
	\label{table:model_PI}
\end{table}

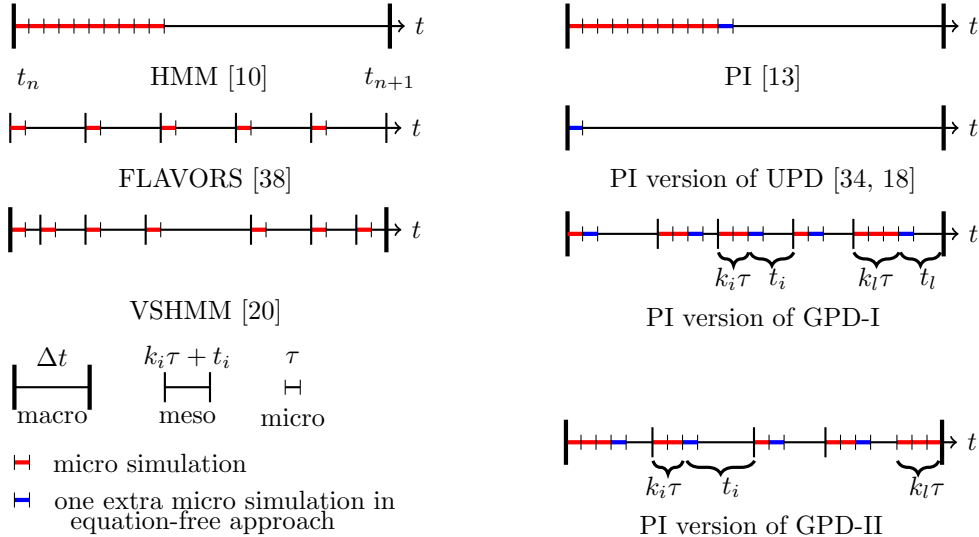
\begin{figure}[h!]
	\begin{tikzpicture}[scale=10]
	% HMM
	\draw[->,thick] (0,0) -- (0.52,0);
	\draw[red,ultra thick] (0,0) -- (0.2,0);
	\draw[ultra thick] (0,-0.03) -- (0,0.03);
	\draw[ultra thick] (0.5,-0.03) -- (0.5,0.03);
	\draw (0.02,-0.01) -- (0.02,0.01);
	\draw (0.04,-0.01) -- (0.04,0.01);
	\draw (0.06,-0.01) -- (0.06,0.01);
	\draw (0.08,-0.01) -- (0.08,0.01);
	\draw (0.10,-0.01) -- (0.10,0.01);
	\draw (0.12,-0.01) -- (0.12,0.01);
	\draw (0.14,-0.01) -- (0.14,0.01);
	\draw (0.16,-0.01) -- (0.16,0.01);
	\draw (0.18,-0.01) -- (0.18,0.01);
	\draw (0.20,-0.01) -- (0.20,0.01);
	\draw (0.26,-0.07) node {HMM \cite{2007_engquist_heterogeneous}};
	\draw[thick] (0.54,0) node {$t$};
	\draw (0.02,-0.07) node {$t_n$};
	\draw (0.5,-0.07) node {$t_{n+1}$};
	\end{tikzpicture}\hspace{0.5cm}
	\begin{tikzpicture}[scale=10]
	%---------------------------
	% PI_Kev
	\draw[->,thick] (0,0) -- (0.52,0);
	\draw[red,ultra thick] (0,0) -- (0.2,0);
	\draw[blue,ultra thick] (0.2,0) -- (0.22,0);
	\draw[ultra thick] (0,-0.03) -- (0,0.03);
	\draw[ultra thick] (0.5,-0.03) -- (0.5,0.03);
	\draw (0.02,-0.01) -- (0.02,0.01);
	\draw (0.04,-0.01) -- (0.04,0.01);
	\draw (0.06,-0.01) -- (0.06,0.01);
	\draw (0.08,-0.01) -- (0.08,0.01);
	\draw (0.10,-0.01) -- (0.10,0.01);
	\draw (0.12,-0.01) -- (0.12,0.01);
	\draw (0.14,-0.01) -- (0.14,0.01);
	\draw (0.16,-0.01) -- (0.16,0.01);
	\draw (0.18,-0.01) -- (0.18,0.01);
	\draw (0.20,-0.01) -- (0.20,0.01);
	\draw (0.22,-0.01) -- (0.22,0.01);
	\draw[thick] (0.26,-0.07) node {PI \cite{2003_gear_projective}};
	\draw[thick] (0.54,0) node {$t$};
	\end{tikzpicture}
	\vspace{0cm}
	\begin{tikzpicture}[scale=10]
	% FLAVORS
	\draw[->,thick] (0,0) -- (0.52,0);
	%\draw[red,ultra thick] (0,0) -- (0.2,0);
	\draw[thick] (0,-0.02) -- (0,0.02);
	\draw[thick] (0.5,-0.02) -- (0.5,0.02);
	\draw[red,ultra thick] (0,0) -- (0.02,0);
	\draw[red,ultra thick] (0.1,0) -- (0.12,0);
	\draw[red,ultra thick] (0.2,0) -- (0.22,0);
	\draw[red,ultra thick] (0.3,0) -- (0.32,0);
	\draw[red,ultra thick] (0.4,0) -- (0.42,0);
	\draw[thick] (0.1,-0.02) -- (0.1,0.02);
	\draw[thick] (0.2,-0.02) -- (0.2,0.02);
	\draw[thick] (0.3,-0.02) -- (0.3,0.02);
	\draw[thick] (0.4,-0.02) -- (0.4,0.02);
	\draw (0.02,-0.01) -- (0.02,0.01);
	\draw (0.12,-0.01) -- (0.12,0.01);
	\draw (0.22,-0.01) -- (0.22,0.01);
	\draw (0.32,-0.01) -- (0.32,0.01);
	\draw (0.42,-0.01) -- (0.42,0.01);
	\draw[thick] (0.54,0) node {$t$};
	\draw[thick] (0.26,-0.07) node {FLAVORS \cite{2010_tao_nonintrusive}};
	\end{tikzpicture}
	\hspace{0.5cm}
	\begin{tikzpicture}[scale=10]
	%---------------------------
	% PI_Samaye
	\draw[->,thick] (0,0) -- (0.52,0);
	\draw[blue,ultra thick] (0,0) -- (0.02,0);
	\draw[ultra thick] (0,-0.03) -- (0,0.03);
	\draw[ultra thick] (0.5,-0.03) -- (0.5,0.03);
	\draw (0.02,-0.01) -- (0.02,0.01);
	\draw[thick] (0.54,0) node {$t$};
	\draw[thick] (0.26,-0.07) node {PI version of UPD \cite{2006_samaey_patch,2009_kevrekidis_equation}};
	\end{tikzpicture}
	\vspace{0cm}
	\begin{tikzpicture}[scale=10]
	% VSHMM
	\draw[->,thick] (0,0) -- (0.52,0);
	%\draw[red,ultra thick] (0,0) -- (0.2,0);
	
	\draw[red,ultra thick] (0,0) -- (0.02,0);
	\draw[red,ultra thick] (0.04,0) -- (0.06,0);
	\draw[red,ultra thick] (0.1,0) -- (0.12,0);
	\draw[red,ultra thick] (0.18,0) -- (0.2,0);
	\draw[red,ultra thick] (0.32,0) -- (0.34,0);
	\draw[red,ultra thick] (0.40,0) -- (0.42,0);
	\draw[red,ultra thick] (0.46,0) -- (0.48,0);
	%\draw[red,ultra thick] (0.48,0) -- (0.5,0);
	\draw[ultra thick] (0,-0.03) -- (0,0.03);
	\draw[ultra thick] (0.5,-0.03) -- (0.5,0.03);
	\draw[thick] (0.04,-0.02) -- (0.04,0.02);
	\draw[thick] (0.1,-0.02) -- (0.1,0.02);
	\draw[thick] (0.18,-0.02) -- (0.18,0.02);
	\draw[thick] (0.32,-0.02) -- (0.32,0.02);
	\draw[thick] (0.4,-0.02) -- (0.4,0.02);
	\draw[thick] (0.46,-0.02) -- (0.46,0.02);
	%	\draw[thick] (0.48,-0.02) -- (0.48,0.02);
	\draw (0.02,-0.01) -- (0.02,0.01);
	\draw (0.06,-0.01) -- (0.06,0.01);
	\draw (0.12,-0.01) -- (0.12,0.01);
	\draw (0.2,-0.01) -- (0.2,0.01);
	\draw (0.34,-0.01) -- (0.34,0.01);
	\draw (0.42,-0.01) -- (0.42,0.01);
	\draw (0.48,-0.01) -- (0.48,0.01);
	\draw[thick] (0.54,0) node {$t$};
	\draw[thick] (0.26,-0.11) node {VSHMM \cite{2013_lee_variable}};
	\end{tikzpicture}\hspace{1.7cm}
	\begin{tikzpicture}[scale=10]
	%---------------------------
	% PI_GPD-I
	
	\draw[->,thick] (0,0) -- (0.52,0);
	\draw[ultra thick] (0,-0.03) -- (0,0.03);
	\draw[ultra thick] (0.5,-0.03) -- (0.5,0.03);
	\draw[red,ultra thick] (0,0) -- (0.02,0);
	\draw[blue,ultra thick] (0.02,0) -- (0.04,0);
	%	\draw[red,ultra thick] (0.1,0) -- (0.12,0);
	\draw[red,ultra thick] (0.12,0) -- (0.14,0);
	\draw[red,ultra thick] (0.14,0) -- (0.16,0);
	\draw[blue,ultra thick] (0.16,0) -- (0.18,0);
	\draw[red,ultra thick] (0.2,0) -- (0.22,0);
	\draw[red,ultra thick] (0.22,0) -- (0.24,0);
	\draw[blue,ultra thick] (0.24,0) -- (0.26,0);
	\draw[red,ultra thick] (0.3,0) -- (0.32,0);
	\draw[blue,ultra thick] (0.32,0) -- (0.34,0);
	\draw[red,ultra thick] (0.38,0) -- (0.4,0);
	\draw[red,ultra thick] (0.4,0) -- (0.42,0);
	\draw[red,ultra thick] (0.42,0) -- (0.44,0);
	\draw[blue,ultra thick] (0.44,0) -- (0.46,0);
	\draw[thick] (0.12,-0.02) -- (0.12,0.02);
	\draw[thick] (0.2,-0.02) -- (0.2,0.02);
	\draw[thick] (0.3,-0.02) -- (0.3,0.02);
	\draw[thick] (0.38,-0.02) -- (0.38,0.02);
	\draw (0.02,-0.01) -- (0.02,0.01);
	\draw (0.04,-0.01) -- (0.04,0.01);
	%	\draw (0.12,-0.01) -- (0.12,0.01);
	\draw (0.14,-0.01) -- (0.14,0.01);
	\draw (0.16,-0.01) -- (0.16,0.01);
	\draw (0.18,-0.01) -- (0.18,0.01);
	\draw (0.22,-0.01) -- (0.22,0.01);
	\draw (0.24,-0.01) -- (0.24,0.01);
	\draw (0.26,-0.01) -- (0.26,0.01);
	\draw (0.32,-0.01) -- (0.32,0.01);
	\draw (0.34,-0.01) -- (0.34,0.01);
	\draw (0.4,-0.01) -- (0.4,0.01);
	\draw (0.42,-0.01) -- (0.42,0.01);
	\draw (0.44,-0.01) -- (0.44,0.01);
	\draw (0.46,-0.01) -- (0.46,0.01);
	%\draw (0.02,0.14) -- (0.02,0.16);
	\draw[thick] (0.54,0) node {$t$};
	\draw[thick] (0.26,-0.11) node {PI version of GPD-I};
	
	\draw [very thick, decorate,decoration={brace,mirror,amplitude=5pt},xshift=5pt,yshift=3pt]
	(0.025,-0.13) -- (0.065,-0.13) node [blue,midway,yshift=0cm] 
	{\footnotesize };
	\draw[thick] (0.22,-0.06) node {$k_i\tau$};
	\draw [very thick, decorate,decoration={brace,mirror,amplitude=5pt},xshift=5pt,yshift=3pt]
	(0.067,-0.13) -- (0.125,-0.13) node [blue,midway,yshift=0cm] 
	{\footnotesize };
	\draw[thick] (0.28,-0.06) node {$t_i$};
	\draw [very thick, decorate,decoration={brace,mirror,amplitude=5pt},xshift=5pt,yshift=3pt]
	(0.205,-0.13) -- (0.265,-0.13) node [blue,midway,yshift=0cm] 
	{\footnotesize };
	\draw[thick] (0.41,-0.06) node {$k_l\tau$};
	\draw [very thick, decorate,decoration={brace,mirror,amplitude=5pt},xshift=5pt,yshift=3pt]
	(0.267,-0.13) -- (0.32,-0.13) node [blue,midway,yshift=0cm] 
	{\footnotesize };
	\draw[thick] (0.48,-0.06) node {$t_l$};
	\end{tikzpicture}
	\vspace{0cm}
	\begin{tikzpicture}[scale=10]
	\draw[thick] (0,0) -- (0.1,0);
	\draw[ultra thick] (0,-0.03) -- (0,0.03);
	\draw[ultra thick] (0.1,-0.03) -- (0.1,0.03);
	\draw[ultra thick] (0.05,0.04) node {$\Delta t$};
	\draw[ultra thick] (0.05,-0.04) node {macro};
	
	%	\draw[ultra thick] (0.12,-0.01) node {\rlap{,}};
	
	\draw[thick] (0.2,0) -- (0.26,0);
	\draw[thick] (0.2,-0.02) -- (0.2,0.02);
	\draw[thick] (0.26,-0.02) -- (0.26,0.02);
	\draw[ultra thick] (0.23,0.04) node {$k_i\tau+t_i$};
	\draw[ultra thick] (0.23,-0.04) node {meso};
	
	%	\draw[ultra thick] (0.22,-0.01) node {\rlap{,}};
	
	\draw[thick] (0.36,0) -- (0.38,0);
	\draw (0.36,-0.01) -- (0.36,0.01);
	\draw (0.38,-0.01) -- (0.38,0.01);
	\draw[ultra thick] (0.37,0.04) node {$\tau$};
	\draw[ultra thick] (0.37,-0.04) node {micro};
	
	%	\draw[ultra thick] (0.22,-0.01) node {\rlap{,}};
	
	\draw[red,ultra thick] (0,-0.1) -- (0.02,-0.1);
	\draw (0,-0.11) -- (0,-0.09);
	\draw (0.02,-0.11) -- (0.02,-0.09);
	\draw[ultra thick] (0.18,-0.1) node {$\text{micro simulation}$};
	
	\draw[blue,ultra thick] (0,-0.15) -- (0.02,-0.15);
	\draw (0,-0.14) -- (0,-0.16);
	\draw (0.02,-0.14) -- (0.02,-0.16);
	\draw[ultra thick] (0.28,-0.15) node {$\text{one extra micro simulation in}$};
	\draw[ultra thick] (0.25,-0.18) node {$\text{equation-free approach}$};
	\end{tikzpicture}\hspace{2cm}
	\begin{tikzpicture}[scale=10]
	% PI_GPD-II	
	\draw[->,thick] (0,0) -- (0.52,0);
	\draw[red,ultra thick] (0,0) -- (0.02,0);
	\draw[red,ultra thick] (0.02,0) -- (0.04,0);
	\draw[red,ultra thick] (0.04,0) -- (0.06,0);
	\draw[blue,ultra thick] (0.06,0) -- (0.08,0);
	\draw[red,ultra thick] (0.115,0) -- (0.135,0);
	\draw[red,ultra thick] (0.135,0) -- (0.155,0);
	\draw[blue,ultra thick] (0.155,0) -- (0.175,0);
	%	\draw[red,ultra thick] (0.23,0) -- (0.25,0);
	\draw[red,ultra thick] (0.25,0) -- (0.27,0);
	\draw[blue,ultra thick] (0.27,0) -- (0.29,0);
	\draw[red,ultra thick] (0.345,0) -- (0.365,0);
	\draw[red,ultra thick] (0.365,0) -- (0.385,0);
	\draw[blue,ultra thick] (0.385,0) -- (0.405,0);
	\draw[red,ultra thick] (0.44,0) -- (0.46,0);
	\draw[red,ultra thick] (0.46,0) -- (0.48,0);
	\draw[red,ultra thick] (0.48,0) -- (0.50,0);
	\draw[ultra thick] (0,-0.03) -- (0,0.03);
	\draw[ultra thick] (0.5,-0.03) -- (0.5,0.03);
	\draw[thick] (0.115,-0.02) -- (0.115,0.02);
	\draw[thick] (0.25,-0.02) -- (0.25,0.02);
	\draw[thick] (0.345,-0.02) -- (0.345,0.02);
	\draw (0.02,-0.01) -- (0.02,0.01);
	\draw (0.04,-0.01) -- (0.04,0.01);
	\draw (0.06,-0.01) -- (0.06,0.01);
	\draw (0.08,-0.01) -- (0.08,0.01);
	\draw (0.135,-0.01) -- (0.135,0.01);
	\draw (0.155,-0.01) -- (0.155,0.01);
	\draw (0.175,-0.01) -- (0.175,0.01);
	\draw (0.25,-0.01) -- (0.25,0.01);
	\draw (0.27,-0.01) -- (0.27,0.01);
	\draw (0.29,-0.01) -- (0.29,0.01);
	\draw (0.365,-0.01) -- (0.365,0.01);
	\draw (0.385,-0.01) -- (0.385,0.01);
	\draw (0.405,-0.01) -- (0.405,0.01);
	\draw (0.44,-0.01) -- (0.44,0.01);
	\draw (0.46,-0.01) -- (0.46,0.01);
	\draw (0.48,-0.01) -- (0.48,0.01);
	%	\draw (0.46,-0.01) -- (0.46,0.01);
	%	\draw (0.02,0.14) -- (0.02,0.16);
	\draw[thick] (0.54,0) node {$t$};
	\draw[thick] (0.26,-0.11) node {PI version of GPD-II};
	
	\draw [very thick, decorate,decoration={brace,mirror,amplitude=5pt},xshift=5pt,yshift=3pt]
	(-0.06,-0.13) -- (-0.02,-0.13) node [blue,midway,yshift=0cm] 
	{\footnotesize };
	\draw[thick] (0.13,-0.06) node {$k_i\tau$};
	\draw [very thick, decorate,decoration={brace,mirror,amplitude=5pt},xshift=5pt,yshift=3pt]
	(-0.015,-0.13) -- (0.075,-0.13) node [blue,midway,yshift=0cm] 
	{\footnotesize };
	\draw[thick] (0.22,-0.06) node {$t_i$};
	\draw [very thick, decorate,decoration={brace,mirror,amplitude=5pt},xshift=5pt,yshift=3pt]
	(0.265,-0.13) -- (0.32,-0.13) node [blue,midway,yshift=0cm] 
	{\footnotesize };
	\draw[thick] (0.48,-0.06) node {$k_l\tau$};
	\end{tikzpicture}
	\caption{Schematics of HMM, FLAVORS, VSHMM, PI, PI version of UPD and PI version of GPD schemes of both types.} \label{fig:PI}
\end{figure}

In figure \ref{fig:PI}, maintaining the same as the stucture as the PI-GPD scheme (except one extra micro simulation), HMM techniques can be applied based on their own macro to micro operator, micro to macro operator, inner \& outer integrators, and bridging techniques. To do large time integration in the HMM approach, neighboring macroscopic informations are utilized, instead of an additional microsimulation (blue segment) in the equation-free approach. HMM \cite{2003_engquist_heterognous,2003_vanden_numerical,2007_engquist_heterogeneous} and it's variants FLAVORS \cite{2010_tao_nonintrusive}, VSHMM \cite{2013_lee_variable}, BA \cite{2009_weinan_general,2015_maclean_note} can also be used in the GPD scheme with a lot of flexibility inside the structure as shown in the following Table \ref{table:model_HMM}.

\begin{table}[h!]
	\begin{center}
		\begin{tabular}{|p{3cm}|p{3.5cm}|p{5cm}|}
			\hline
			%			%\cline{2-3}
			%			&\multicolumn{2}{|c|}{Distribution of 24 number of gap-}\\
			%			&\multicolumn{2}{|c|}{tooth timesteppers in GPD scheme of}\\
			%\hline
			\textbf{Proposed scheme}&\textbf{Values of $l$ and / or $k$}&\textbf{The proposed scheme will be equivalent to}\\
			\hline
			\multirow{2}{*}{Uniform PI-GPD-I}&$l=1$&
			HMM \cite{2007_engquist_heterogeneous}\\\cline{2-3}
			&$k_i=1$&FLAVORS \cite{2010_tao_nonintrusive}\\
			\hline
			Non-uniform PI-GPD-I&$k_i=1$&VSHMM \cite{2013_lee_variable}\\
			\hline
		\end{tabular}
	\end{center}
	\caption{The PI-GPD scheme is a generalized version of HMM, FLAVORS and VSHMM schemes.}
	\label{table:model_HMM}
\end{table}

By putting $k_i=1$, where $i=1, 2,..., l$ for some finite $l$ in the uniform GPD scheme of type-I (except one extra micro simulation in each meso time step) with replacing the stiffness parameter ($\epsilon$) by a higher value, one can obtain the boosting algorithm (BA) in heterogeneous multiscale methods (HMM) as proposed by Maclean \cite{2015_maclean_note}. Using this boosting technique a large computational advantage can be obtained.

%\textcolor{blue}{Many methods are similar in terms of structure for particular choices of $l$ and $k$ such as:}

\section{Results and discussion}

In order to verify the results, an unsteady linear reaction-diffusion problem is considered as follows: 

\begin{equation}\label{eqn:22}
\begin{cases}
\text{PDE:}\hspace{0.25cm} & u_t=u_{xx}+u,\hspace{0.25cm} \forall \hspace{0.25cm}  (x,t)\in(0,1)\times(0,1],\\
\text{IC:}\hspace{0.25cm} & u(x,0)=\sin(\pi x),\hspace{0.25cm} \forall\hspace{0.25cm} x\in(0,1),\\
\text{BCs:}\hspace{0.25cm} & u(0,t)=u(1,t)=0, \hspace{0.25cm} \forall\hspace{0.25cm} t\in[0,1].
\end{cases}
\end{equation}
It is chosen to be a microscopic problem. For a macroscopic discretization, $\Delta x=5e{-2}$ and $\Delta t$=1$e{-3}$ are considered. $Nt$ is the numbers of long time steps of size $\Delta t$ such that $Nt\times \Delta t=\text{final time}$ (=1). To construct the patches centred at the macroscopic grid points, $h=2e{-2}$, $\tau=1e{-5}$ are chosen. The Neumann boundary conditions as given in \eqref{eqn:7} are considered in the microscopic simulation, which is physically relevant with respect to the model problem. The equation \eqref{eqn:8} is used as a lifting operator (micro initial condition) to evolve the microscopic problem \eqref{eqn:22} for short time $\tau$ in the GTT. At time $t=0$, GTT is started with a constant error due to the averaging of the initial condition inside the patches. Here, in the microscopic simulation, one can use a finite difference discretization to evolve the microscopic problem for time $\tau$ inside the patches. For the simulation inside the patches, FTCS scheme is implemented in this study to find all results. To restrict the microscopic solution to the macroscopic solution inside the patches, Simpson's 1/3 rule is used as a numerical quadrature formula. In the UPD scheme, the total $k=24$ number of GTTs are implemented at the starting of each long time step $\Delta t$, and to do a long time integration, the forward Euler method is used as an outer integrator to march forward for time $7.6e{-4}$. In the GPD scheme of both types, total 24 number of GTTs are distributed in all possible ways in a long time step $\Delta t$ to find the solution of the system \eqref{eqn:22}. 

\textbf{Uniform case:} The results are discussed based on the uniform distribution of both GTTs and extrapolation times in each time step $\Delta t$. GPD solutions of both types of the microscopic problem \eqref{eqn:22} are compared with the analytical solution of the effective macroscopic equation of \eqref{eqn:22}. That analytical solution is 

\begin{equation}\label{eqn:23}
U_a(x,t)=\sin(\pi x) e^{(1-\pi ^2)t}. 
\end{equation}

\begin{figure}[h!]
	\centering
	\includegraphics[width=12cm,height=7cm]{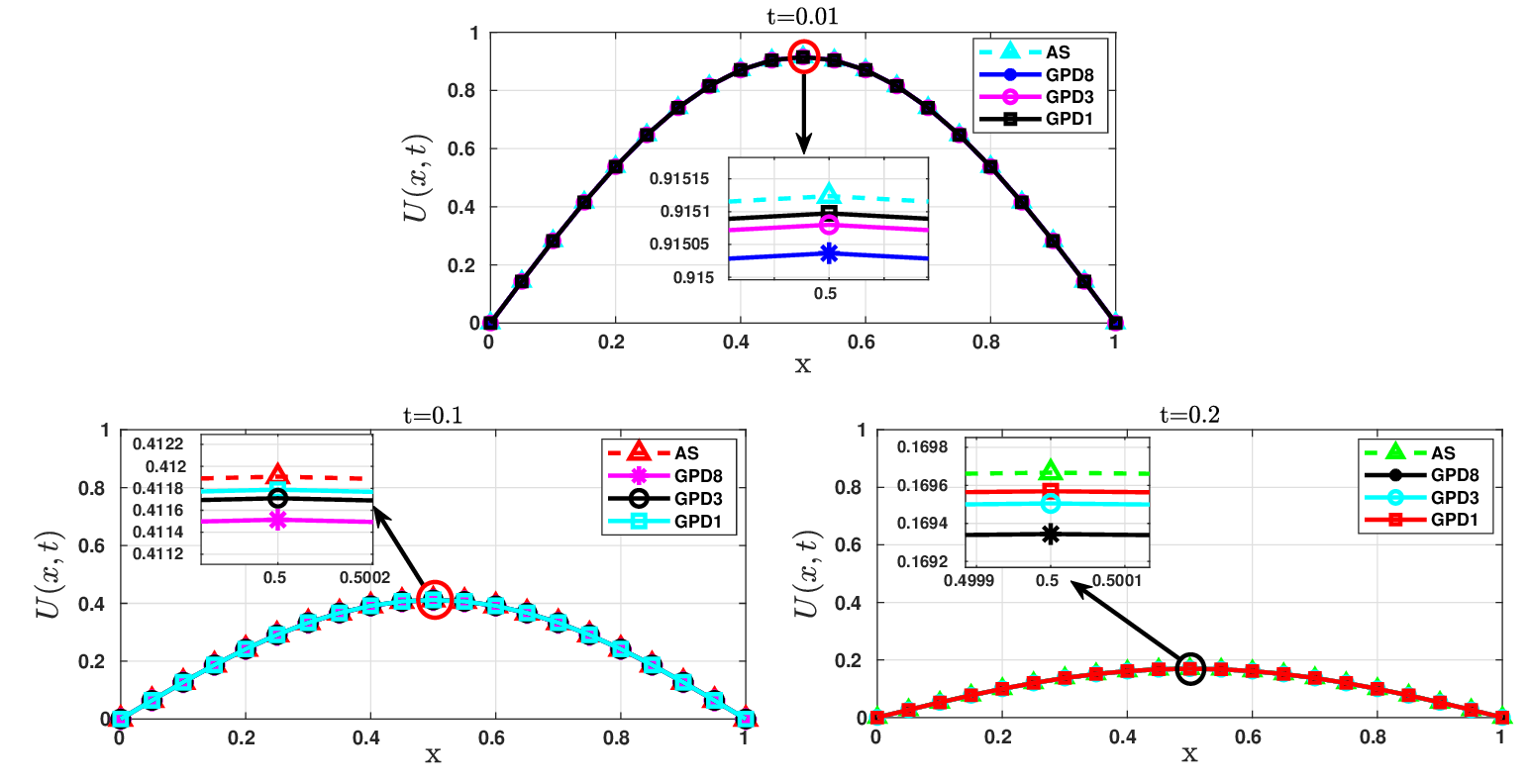}          \caption{Comparison between the analytical solution and the GPD solutions of type-I of the effective macroscopic equation of \eqref{eqn:22} at times $t=0.01, t=0.1$ and $t=0.2$ are shown in this figure. }
	\label{fig:5_XvsU_Type1_Old}
\end{figure}

%\pagebreak

\begin{figure}[h!]
	\centering
	\includegraphics[width=12cm,height=7cm]{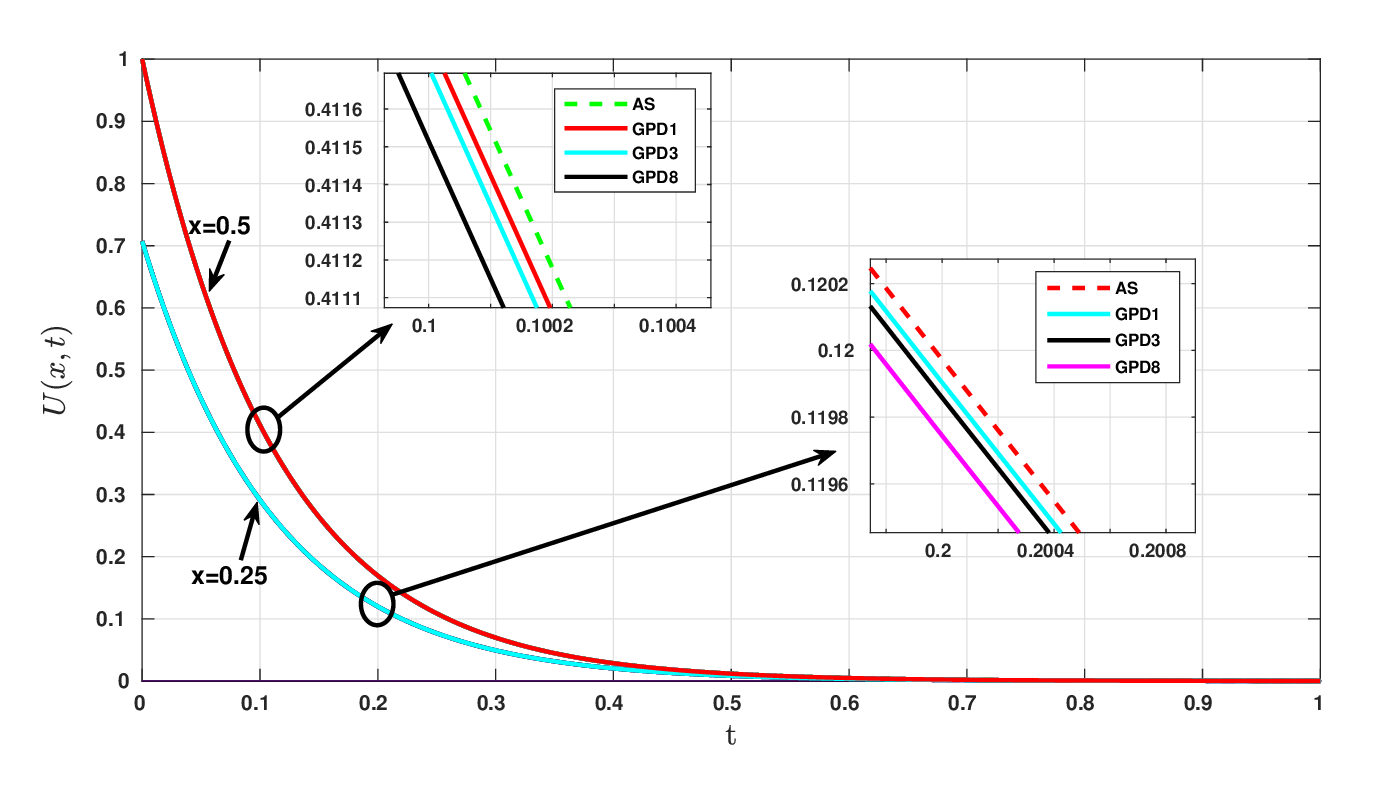}          \caption{Comparison between the analytical solution and the GPD solutions of type-I of the effective macroscopic equation of \eqref{eqn:22} at $x=0.25$ and $x=0.5$ are shown in this figure. }
	\label{fig:6_TvsU_Type1_Old}
\end{figure}

In order to compare with the analytical solution, the patch is discretized at the microscopic level. We take $dx=4.4444e-5$ and $dt=9.5238e-10$. Total 24 number of GTTs are distributed in three different ways in the GPD scheme, such as \{8, 8, 8\}, \{3, 3,..., 3\} and \{1,1,...,1\} in each time step $\Delta t$, and the corresponding GPD solutions of type-I are denoted by GPD8, GPD3 and GPD1 respectively. The uniform distribution of GTTs and uniform extrapolation time steps are considered to find the GPD solutions of type-I. It can be observed from the above figures that when the GTTs are distributed over the time step $\Delta t$, extrapolation time decreases, numbers of GTTs and number of spatial bridgings are increased at different time locations in $\Delta t$ time step. Consequently, better accuracy is achieved. As the number of distributions increases, the GPD solutions (GPD8, GPD3, and GPD1) become more closer to the analytical solution. GTTs of GPD3 are more sparsely distributed than those of GPD8, effectively GPD3 gives better accuracy compared to the GPD8. Similarly, the GTTs of GPD1 are even more sparsely distributed, hence GPD1 has a better accuracy than both GPD8 and GPD3. Similar kinds of accuracy is achieved against the analytical solution for the GPD scheme of type-II. Figures \ref{fig:5_XvsU_Type1_Old} \& \ref{fig:6_TvsU_Type1_Old} show a comparison between the GPD solutions of type-I of the problem \eqref{eqn:22} and the analytical solution of the effective macroscopic problem of \eqref{eqn:22}.

\begin{table}[h!]
	\begin{center}
		\scalebox{0.7}{
			\begin{tabular}{ |p{2cm}||p{1.5cm}|p{1.5cm}|p{1.5cm}|p{1.6cm}||p{1.5cm}|p{1.5cm}|p{1.5cm}|p{1.6cm}|}
				\hline
				&\multicolumn{4}{|c||}{GPD Scheme of Type-I} & \multicolumn{4}{c|}{GPD Scheme of Type-II} \\
				\hline
				Distribution of GTTs in $\Delta t$ Time Step&No of GTTs &No of Extrapolations&Size of Extrapolations&Maximum \% Error&No of GTTs & No of Extrapolations&Size of Extrapolations&Maximum \% Error\\
				\hline
				24 (UPD)&25&1&$\Delta t-24\tau$&4.3755&&&&\\
				\hline
				\{12, 12\}&26&2&$\frac{\Delta t-24\tau}{2}$&2.0229&25&1&$\Delta t-24\tau$&4.3755\\
				\{8, 8, 8\}&27&3&$\frac{\Delta t-24\tau}{3}$&1.2304&26&2&$\frac{\Delta t-24\tau}{2}$&2.0229\\
				\{6, 6, 6, 6\}&28&4&$\frac{\Delta t-24\tau}{4}$&0.8326&27&3&$\frac{\Delta t-24\tau}{3}$&1.2304\\
				\{5, 5, 5, 5, 4\}&29&5&$\frac{\Delta t-24\tau}{5}$&0.5934&28&4&$\frac{\Delta t-24\tau}{4}$&0.8326\\
				\{4,...4\}&30&6&$\frac{\Delta t-24\tau}{6}$&0.4338&29&5&$\frac{\Delta t-24\tau}{5}$&0.5934\\
				\{3,...,3\}&32&8&$\frac{\Delta t-24\tau}{8}$&0.2340&31&7&$\frac{\Delta t-24\tau}{7}$&0.3196\\
				\{2,...,2\}&36&12&$\frac{\Delta t-24\tau}{12}$&0.0339&35&11&$\frac{\Delta t-24\tau}{11}$&0.0703\\
				\{1,...,1\}&48&24&$\frac{\Delta t-24\tau}{24}$&0.1664&47&23&$\frac{\Delta t-24\tau}{23}$&0.1577\\
				\hline
		\end{tabular}}
	\end{center}
	\caption{Total 24 number of GTTs are distributed in few ways over time step $\Delta t$ with $Nt=600$ to compare the accuracy in both GPD schemes.}
	\label{table:1}
\end{table}	

%\begin{table}
%	\begin{center}
%		\scalebox{0.9}{
%			\begin{tabular}{r|lll}
%				\multicolumn{1}{r}{}
%				& \multicolumn{1}{l}{Heading 1}
%				& \multicolumn{1}{l}{Heading 2}
%				& \multicolumn{1}{l}{Heading 3} \\ \cline{2-4}
%				Row 1 & Cell 1,1 & Cell 1,2 & Cell 1,3 \\
%				Row 2 & Cell 2,1 & Cell 2,2 & Cell 2,3
%		\end{tabular}}
%	\end{center}
%\end{table}

In order to do a long time integration as an outer integrator, the forward Euler method (equation \eqref{eqn:12}) is implemented. In Table \ref{table:1}, a comparison is made between the GPD scheme of types I and II for different distributions of GTTs, where $Nt=600$ and all other parameters are kept fixed. It is observed that for fixed $Nt=600$, the GPD scheme of type-I with $k=\{12, 12\}$ has a better accuracy compared to the UPD scheme with $k=24$, because one extra GTT and one extra spatial bridging are computed in the GPD scheme of type-I and the extrapolation time of GPD scheme is half of the extrapolation time of UPD scheme in time step $\Delta t$. Whereas the GPD scheme of type-II with $k=\{12, 12\}$ and UPD scheme with $k=24$ have the same number of GTTs, same number of spatial bridgings, the same number and size of extrapolations in time step $\Delta t$, so both solutions have the similar accuracy. If the number of GTTs, number of spatial bridgings, number \& size of extrapolations are the same, similar accuracy is achieved in both the GPD schemes as shown in Table \ref{table:1} (first six rows), otherwise, they will differ, confirming that the same proposal was made by Maclean et al. \cite{2015_maclean_convergence}. As for example, the accuracy of GPD scheme of type-I with \{4, 4,..., 4\} and GPD scheme of type-II with \{3, 3,..., 3\} are different. GPD schemes of type-I and type-II have different number of extrapolations as 6 and 7, respectively, and hence different extrapolation times of $\frac{Dt-24\tau}{6}$ and $\frac{Dt-24\tau}{7}$ respectively in each $\Delta t$ time step. As the GPD scheme of type-II has a smaller extrapolation time step as well as one extra computation of both GTT and spatial bridging compared to the GPD scheme of type-I, so GPD scheme of type-II has better accuracy. Whatever we expected in the GPD schemes, we have achieved the same. One interesting observation from Table \ref{table:1} is that when the distribution of the GTTs is refined, a better accuracy is observed in the solution. However, this procedure does not work in the last row because the micro-state does not reach its quasi-equilibrium state for the distribution $\{1,...,1\}$. For $\{2,...,2\}$ distribution of GTTs, the micro simulation reaches to it's quasi-equilibrium state and the relaxation time can be determined from the distribution $\{2,...,2\}$ such that $\tau_\epsilon=2\tau=2e-5.$ Hence one needs to choose a coarser distribution of $\{1,...,1\}$ GTT to allow the micro-state to attain its quasi-equilibrium state for this problem.

\begin{table}[h!]
	\begin{center}
		\begin{tabular}{|c|c||c|c|c||c|c|}
			\hline
			\multicolumn{2}{|c||}{\textbf{\large UPD Scheme}} & \multicolumn{3}{c||}{\textbf{\large GPD Scheme of Type-I}} &  \multicolumn{2}{c|}{\textbf{\large Ratio}}\\
			\hline
			\textbf{Nt}&\textbf{Maximum}&\textbf{Nt}&\textbf{Distribution}&\textbf{Maximum}& \textbf{$\text{\%E}_\text{UPD}$/} &\textbf{$\text{T}_\text{UPD}$/}\\
			&\textbf{\% Error}&&\textbf{of GTTs}&\textbf{\% Error}&\textbf{$\text{\%E}_\text{GPD}$}&\textbf{$\text{T}_\text{GPD}$}\\
			&(\textbf{$\text{\%E}_\text{UPD}$})&&&(\textbf{$\text{\%E}_\text{GPD}$})&  & \\
			\hline
			1000&1.8995&800&\{12,12\}&1.2363&1.5364&0.9992\\
			%&&700&\{12, 12\}&1.5722&1.2082&1.1002\\
			&&100&\{1, 1,...,1\}&1.1855&1.6023&1.6051\\
			\hline
			%900&2.3075&600&\{12, 12\}&2.0229&1.1407&1.1578\\
			900&2.3075&100&\{1, 1,..., 1\}&1.1855&1.9464&1.4468\\
			\hline
			%800&2.8216&600&\{12, 12\}&2.0229&1.3948&1.0332\\
			%&&500&\{12, 12\}&2.6562&1.0623&1.2214\\
			800&2.8216&100&\{2, 2,..., 2\}&2.7260&1.0351&2.2159\\
			&&100&\{1, 1,..., 1\}&1.1855&2.3801&1.2910\\
			\hline
			700&3.4864&600&\{12, 12\}&2.0229&1.7235&0.9925\\
			%&&500&\{12, 12\}&2.6562&1.3126&1.1713\\
			&&300&\{8, 8, 8\}&3.3569&1.0386&1.9218\\
			%&&200&\{5, 5, 5, 5, 4\}&3.1610&1.1029&1.2749\\
			&&100&\{2, 2,..., 2\}&2.7260&1.2789&2.1285\\
			&&100&\{1, 1,..., 1\}&1.1855&2.9409&1.2401\\
			\hline
			600&4.3755&500&\{12, 12\}&2.6562&1.6473&1.0038\\
			%&&400&\{12, 12\}&3.6069&1.2131&1.2543\\
			%&&300&\{8, 8, 8\}&3.3569&1.3034&1.6470\\
			&&200&\{6, 6, 6, 6\}&4.0304&1.0856&2.2371\\
			%&&200&\{5, 5, 5, 5, 4\}&3.1610&1.3842&1.0927\\
			&&100&\{3, 3,..., 3\}&4.2498&1.0296&3.8898\\
			&&100&\{2, 2,..., 2\}&2.7260&1.6051&1.8243\\
			&&100&\{1, 1,..., 1\}&1.1855&3.6908&1.0629\\
			\hline
		\end{tabular}
	\end{center}
	\caption{Comparison between UPD scheme and GPD scheme of Type-I. }
	\label{table:2}
\end{table}

%\pagebreak
From Tables \ref{table:2} \& \ref{table:3}, it can be observed that the percentage errors and computational time in the GPD schemes are improved compared to those of the UPD scheme. $\text{\%E}_\text{UPD}$ and $\text{\%E}_\text{GPD}$ are the maximum percentage errors in the UPD scheme and GPD scheme respectively in the whole domain and over the full time interval [0, 1]. $\text{T}_\text{UPD}$ and $\text{T}_\text{GPD}$ are the computational times in UPD scheme and GPD scheme respectively. For an example from the Table \ref{table:2}, UPD scheme with $Nt=1000$ and the GPD scheme of type-I with $Nt=100$ along with distribution of GTTs \{1, 1,..., 1\} are considered. GPD scheme of type-I achieves $\frac{5}{8}$ of the percentage error in less than $\frac{5}{8}$ of the time required by the UPD scheme.
%The maximum percentage error in the UPD solution is approximately 1.60 times the maximum percentage error in the GPD solution, and the GPD scheme of type-I has approximately 1.61 times faster convergence rate than the UPD scheme. 
For the UPD scheme, a total of 25 number of GTTs and same number of spatial bridgings are computed in each time step $\Delta t$ around each interior macroscopic grid point $x_j$, $ j=1, 2,..., N-1$. That implies a total of $1000\times 25=25000$ GTTs and same number of spatial bridgings are computed with respect to each interior macroscopic grid point $x_j$ over the entire time interval $[0,1]$ in UPD scheme. Whereas in the GPD scheme of type-I, a total of 48 numbers of GTTs and same number of spatial bridgings are computed over $\Delta t$ time step, that means total $100\times 48=$ 4800 numbers of GTTs as well as same number of spatial bridgings are computed over the entire time interval $[0,1]$ for each interior macroscopic grid point $x_j$. GPD scheme uses less GTTs and spatial bridgings than UPD scheme, resulting in faster computation. Furthermore, instead of a single extrapolation time of length $\Delta t-24\tau$ in the UPD scheme, 24 extrapolations each of size $\frac{\Delta t-24\tau}{24}$ in the GPD scheme of type-I are computed in the long time step $\Delta t$. Hence, a total of $1000\times 1$= 1000 number of extrapolations are implemented in UPD scheme, while $100\times 24$= 2400 number of extrapolations are implemented in the GPD scheme of type-I over the entire time interval $\left[0,1\right]$ for each interior macroscopic grid point $x_j$. Due to the less number of GTTs and spatial bridgings in the GPD scheme, it becomes faster than the UPD scheme. Additionally, due to shorter extrapolation times in the GPD scheme, it reduces the error in the extrapolation. These combined effects make the GPD scheme of type-I more accurate and result in a faster convergence rate than the UPD scheme.

\begin{table}[h!]
	\begin{center}
		\begin{tabular}{|c|c||c|c|c||c|c|}
			\hline
			\multicolumn{2}{|c||}{\textbf{\large UPD Scheme}} & \multicolumn{3}{c||}{\textbf{\large GPD Scheme of Type-II}} &  \multicolumn{2}{c|}{\textbf{\large Ratio}}\\
			\hline
			\textbf{Nt}&\textbf{Maximum}&\textbf{Nt}&\textbf{Distribution}&\textbf{Maximum}& \textbf{$\text{\%E}_\text{UPD}$/} &\textbf{$\text{T}_\text{UPD}$/}\\
			&\textbf{\% Error}&&\textbf{of GTTs}&\textbf{\% Error}&\textbf{$\text{\%E}_\text{GPD}$}&\textbf{$\text{T}_\text{GPD}$}\\
			&(\textbf{$\text{\%E}_\text{UPD}$})&&&(\textbf{$\text{\%E}_\text{GPD}$})&  & \\
			\hline
			%900&2.3075&900&\{12, 12\}&2.3075&1&1.0058\\
			%&&200&\{3, 3,...,3\}&2.1607&1.0679&0.9981\\
			900&2.3075&100&\{1, 1,...,1\}&1.2528&1.8419&2.0417\\
			\hline
			800&2.8216&500&\{8, 8, 8\}&2.6562&1.0623&1.4015\\
			%&&300&\{5, 5, 5, 5, 4\}&2.4342&1.1591&1.0501\\
			&&100&\{1, 1,..., 1\}&1.2528&2.2522&2.0423\\
			\hline
			700&3.4864&700&\{12, 12\}&3.4864&1&1.1046\\
			%&&500&\{8, 8, 8\}&2.6562&1.3126&1.2002\\
			&&300&\{6, 6, 6, 6\}&3.3569&1.0386&1.8501\\
			%&&200&\{4, 4,..., 4\}&3.1610&1.1029&1.3292\\
			&&100&\{2, 2,..., 2\}&3.0043&1.1605&4.1048\\
			&&100&\{1, 1,..., 1\}&1.2528&2.7829&1.7489\\
			\hline
			600&4.3755&600&\{12, 12\}&4.3755&1&1.1382\\
			%&&500&\{8, 8, 8\}&2.6562&1.6473&1.0290\\
			%&&400&\{8, 8, 8\}&3.6069&1.2131&1.4118\\
			%&&300&\{6, 6, 6, 6\}&3.3569&1.3034&1.5862\\
			&&200&\{5, 5, 5, 5, 4\}&4.0304&1.0856&2.3178\\
			%&&200&\{4, 4,..., 4\}&3.1610&1.3842&1.1396\\
			&&100&\{2, 2,..., 2\}&3.0043&1.4564&3.5194\\
			&&100&\{1, 1,..., 1\}&1.2528&3.4926&1.4995\\
			\hline
		\end{tabular}
	\end{center}
	\caption{Comparison between UPD scheme and GPD scheme of Type-II.}
	\label{table:3}
\end{table}

In a similar fashion as the Table \ref{table:2}, the comparison between the UPD scheme and GPD scheme of type-II is shown in Table \ref{table:3}. As it can be seen, the GPD scheme gives better accuracy in less computational time than the UPD scheme. For example, the UPD scheme with $Nt=700$ (or $Nt=600$) and the GPD scheme of type-II with the same $Nt$ along with distribution of GTTs $\{12, 12\}$ have obviously same solution, because both have same number of GTTs, same number of spatial bridgings, same number and size of extrapolations in time step $\Delta t$, but GPD scheme of type-II is $14\%$ (or $10\%$) faster than the UPD scheme. 

Another example, UPD scheme with $Nt=800$ and GPD scheme with $Nt=100$ along with distribution of GTTs $\{1, 1,..., 1\}$ are considered here. 
GPD scheme of type-II achieves less than $\frac{1}{2}$ of the percentage error in less than $\frac{1}{2}$ of the time of the UPD scheme.
In the UPD scheme, a total of 25 number of GTTs as well as same number of spatial bridgings are computed in each time step $\Delta t$ around each interior macroscopic grid point $x_j, \forall j=1, 2,..., N-1$. That means the total $800\times 25= 20000$ number of GTTs as well as same number of spatial bridgings are computed over the entire time interval $[0,1]$ with respect to each interior macroscopic grid point $x_j$. However, under GPD scheme of type-II with $Nt=100$ and distribution \{1, 1,..., 1\}, 47 numbers of GTTs and same number of spatial bridgings are computed over each $\Delta t$ time step, which implies $100\times 47=$ 4700 GTTs as well as same number of spatial bridgings have been computed during the entire time interval [0,1] corresponding to each interior macroscopic grid point $x_j$. Due to less number of micro simulations and spatial bridgings in GPD scheme of type-II than the UPD scheme, the GPD scheme of type-II has faster computation than the UPD scheme. Furthermore, GPD scheme of type-II computes 23 number of extrapolations each of size $\frac{\Delta t-24\tau}{23}$ in long time step $\Delta t$ instead of single extrapolation of time $\Delta t-24\tau$ in UPD scheme. Hence a total of $800\times 1=$ 800 extrapolations are implemented in the UPD scheme, while $100\times 23=$ 2300 are given in the GPD scheme of type-II over the whole time interval [0,1] for each interior macroscopic grid point $x_j$. Due to the less GTTs and spatial bridgings in the GPD scheme, it becomes faster than the UPD scheme. Additionally, due to shorter extrapolation times in the GPD scheme, it reduces the error in the extrapolation. These combined effects make the GPD scheme of type-II more accurate and faster convergence than the UPD scheme.

\begin{figure}[h!]
	\centering
	\includegraphics[width=6cm,height=4cm]{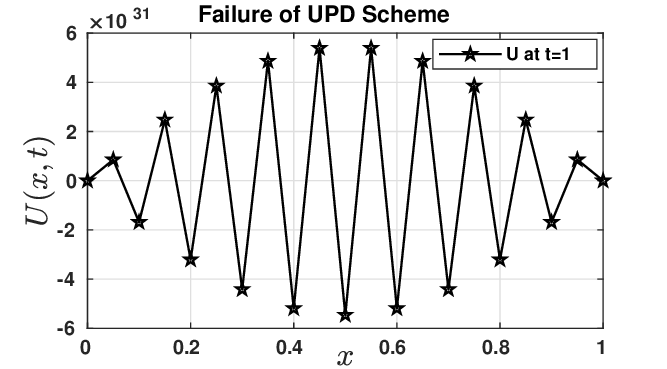} 
	\includegraphics[width=6cm,height=4cm]{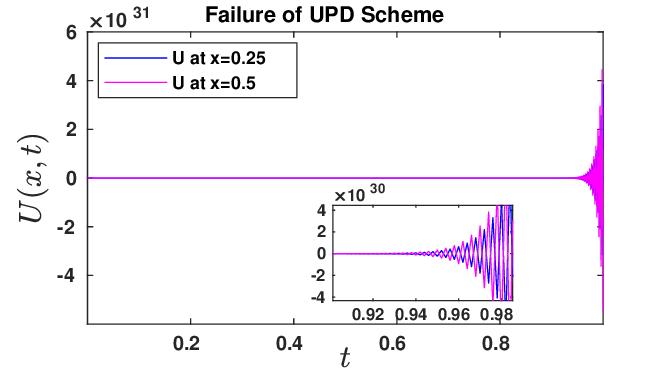}
	\caption{Failure of UPD with $Nt=500$.}
	\label{fig:Failure_1}
\end{figure}

UPD scheme has certain limitations. One of those is the UPD solution of the problem \eqref{eqn:22} diverges for $Nt\le 560$ (approximately) even though the microsimulations converges. That means the UPD scheme converges for macroscopic time step size belongs to $\left(0,\frac{1}{560}\right)$= $\left(0,0.00176\right)$ and diverges for $\left[0.00176, \infty \right)$. That means if the extrapolation time becomes greater or equal than 0.001546 (=0.00176 - 24$\tau$) then the solution diverges in UPD scheme. GPD scheme of both types can resolve such issues of UPD scheme with less error and less computational time. This is an advantage of the GPD scheme over UPD scheme.

%%%%%%%%%%%%%%%%%%%%%%%%%%%%%%%%%%%%%%%%%%%%%%%%%%%%%%%

In figure- \ref{fig:Failure_1}, $Nt=500$ is considered to investigate the effective macroscopic solution of \eqref{eqn:22} using UPD scheme, and all others parameters are kept same as previously declared. Here the extrapolation time is 0.00176 ($> 0.001546$). The UPD solution of \eqref{eqn:22} diverges. In figures \ref{fig:Success_1} \& \ref{fig:Success_2}, the GPD scheme of type-I and type-II are considered with the distribution of 24 number of GTTs as $\{12, 12\}$ and $\{8, 8, 8\}$ respectively and all the other parameters are kept same as previously declared. In both of the GPD schemes the number of extrapolations is two and the extrapolation time is $\frac{0.00176}{2}=0.00088 < 0.001546$ in each scheme, hence the GPD solutions of both types converge. GPD scheme of both types are mostly capable of handling this kind of failures of UPD scheme.

\begin{minipage}{0.39\linewidth}
	\centering
	\includegraphics[width=4.5cm,height=2.5cm]{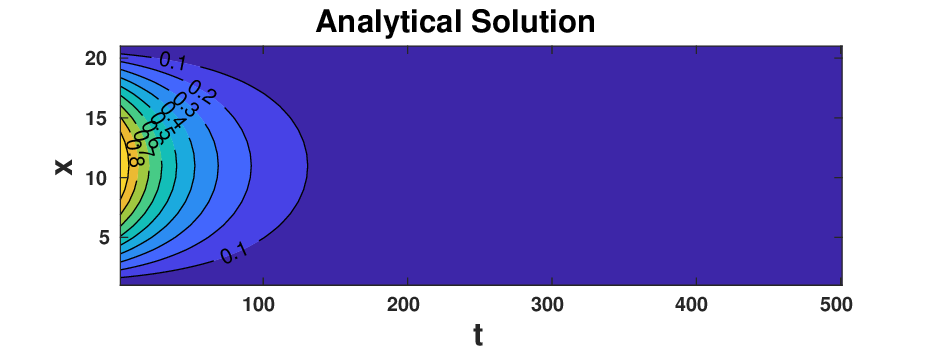}\\
	\includegraphics[width=4.5cm,height=2.5cm]{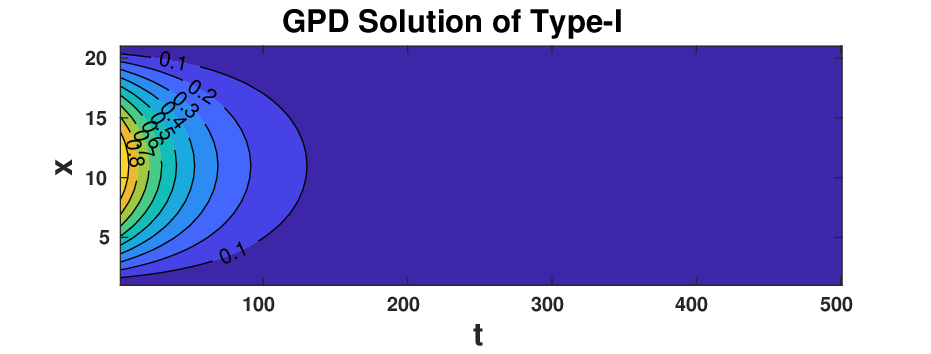}\\
	\includegraphics[width=4.5cm,height=2.5cm]{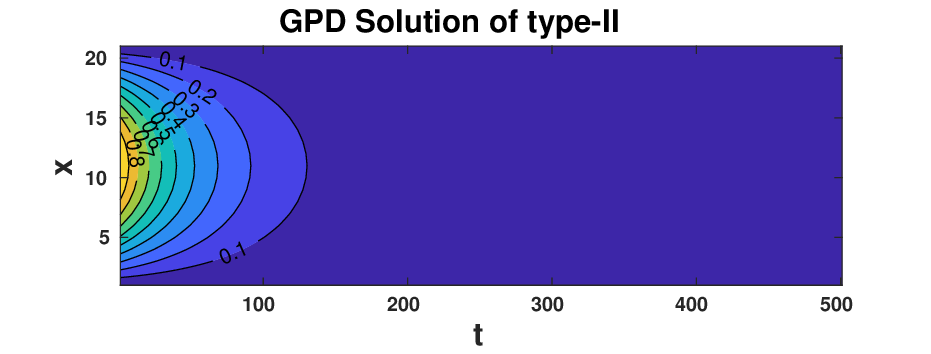}
	\captionof{figure}{Contours of the analytical solution, GPD solution of type-I and GPD solution of type-II for problem \eqref{eqn:22} with $Nt=500$ and the distribution of GTTs in GPD schemes are $\{12, 12\}$ and $\{8, 8, 8\}$.}
	\label{fig:Success_1}
\end{minipage}
\begin{minipage}{0.49\linewidth}
	\centering
	\includegraphics[width=5cm,height=4cm]{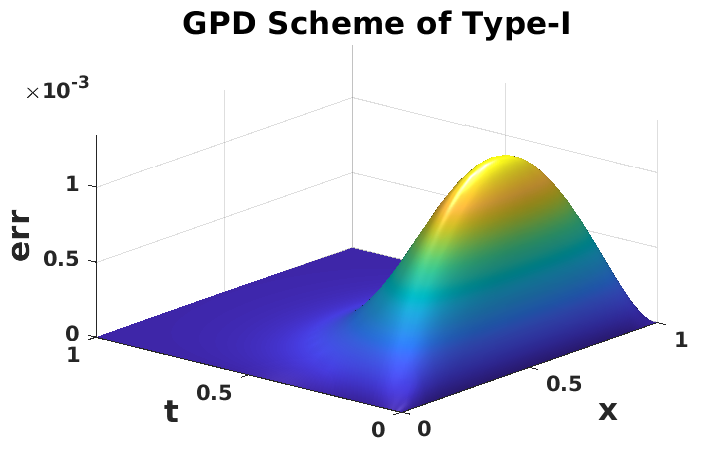}\\
	\includegraphics[width=5cm,height=4cm]{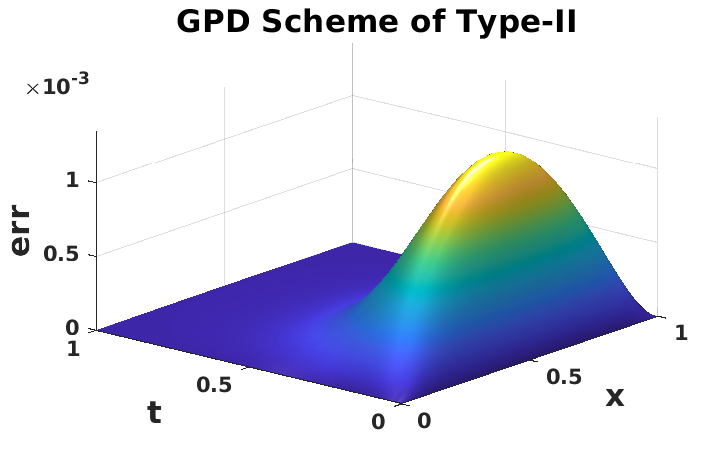}
	\captionof{figure}{Surface plots of errors for problem \eqref{eqn:22}.}
	\label{fig:Success_2}
\end{minipage}

%\pagebreak

\begin{table}[h!]
	\begin{center}
		\begin{tabular}{|p{0.6cm}|p{4cm}|p{4cm}|}
			\hline
			%			%\cline{2-3}
			%			&\multicolumn{2}{|c|}{Distribution of 24 number of gap-}\\
			%			&\multicolumn{2}{|c|}{tooth timesteppers in GPD scheme of}\\
			%\hline
			Nt&Type-I&Type-II\\
			\hline
			600&\{24\}&\{24\}\\
			500&\{12, 12\}&\{8, 8, 8\}\\
			400&\{12, 12\}&\{8, 8, 8\}\\
			300&\{8, 8, 8\}&\{6, 6,..., 6\}\\
			200&\{6, 6,..., 6\}&\{5, 5, 5, 5, 4\}\\
			100&\{3, 3,..., 3\}&\{2, 2,..., 2\}\\
			\hline
		\end{tabular}
	\end{center}
	\caption{Distribution of GTTs in GPD scheme mentioned against a particular $Nt$ or finer than that distribution ensure convergence of solution.}
	\label{table:4}
\end{table}

GPD scheme of both types also have some restrictions to choose appropriate extrapolation time to converge the solution. Due to long extrapolation time if the GPD solution diverges, the GTTs are distributed into more groups within time $\Delta t$. In Table \ref{table:4}, it is shown that all such distributions of total 24 number of GTTs and it's finer distribution that ensures the convergence of the GPD scheme with respect to $Nt$. For an example, when $Nt=100$, GPD scheme of type-I converges only for the distributions $\{3, 3,..., 3\}$, $\{2, 2,..., 2\}$ \& $\{1, 1,..., 1\}$ and GPD scheme of type-II converges only for the distributions $\{2, 2,..., 2\}$ \& $\{1, 1,..., 1\}$.

%%%%%%%%%%%%%%%%%%%%%%%%%%%%%%%%%%%%

\begin{minipage}{0.39\linewidth}
	\centering
	\includegraphics[width=4.5cm,height=2.5cm]{AS_Contour_New.eps}\\
	\includegraphics[width=4.5cm,height=2.5cm]{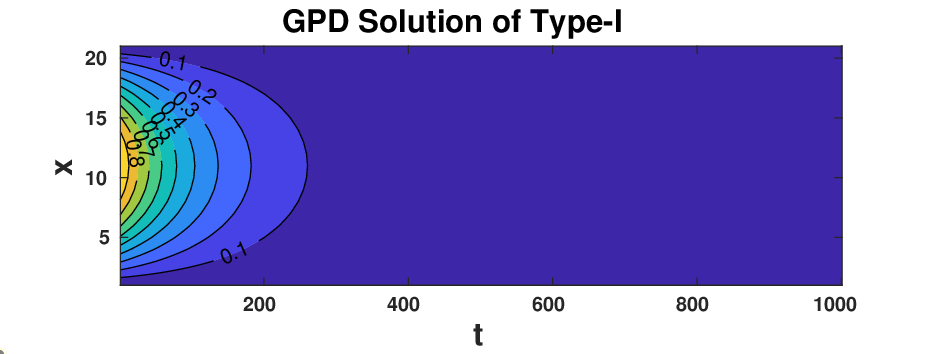}
	\captionof{figure}{Contours of the analytical solution and GPD solution of type-I for problem \eqref{eqn:22} with non-uniform distribution of GTTs and extrapolation times.}
	\label{fig:Contour_T1_NonUni}
\end{minipage}
\begin{minipage}{0.49\linewidth}
	\centering
	\includegraphics[width=5cm,height=4cm]{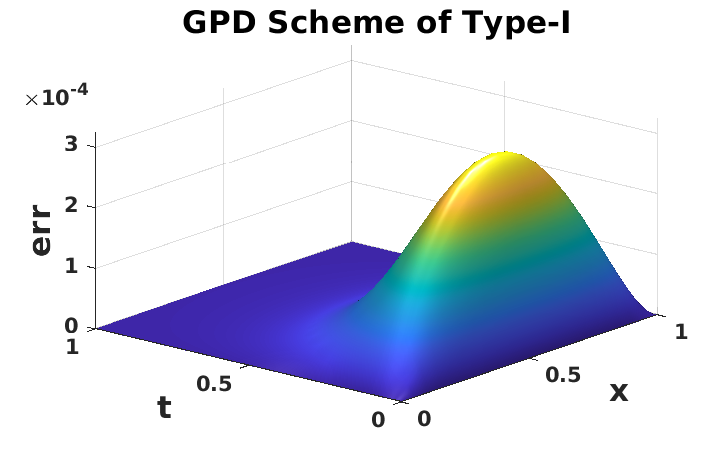}
	\captionof{figure}{Surface plots of errors for problem \eqref{eqn:22}.}
	\label{fig:Contour_T1_NonUni_Error}
\end{minipage}

\vspace{0.5cm}

\textbf{Non-uniform case:} So far only uniform distribution on GTTs and uniform extrapolations have been used. In other words, uniform meso time steps are considered. However, one may need to use variable numbers of GTTs as well as variable extrapolation times (i.e., variable meso time steps) on $\Delta t$ based on the nature of the problem. For an example, chemical plants release waste materials into a river or in the atmosphere over a period of time are common occurrences \cite{2000_mazumder_contaminant}. Such kind of releases are non-uniform in time and mostly periodic. So, for such kind of problems, uniform schemes are not expected to produce desired result and non-uniform schemes are more useful. In such periodic problem, $\Delta t$ can be considered as a period. Like periodic problems, there is no fixed rule for a non-periodic problem to choose $\Delta t$, so $\Delta t$ should be chosen carefully, looking at the nature of the problem so that the high gradient or oscillations zones receives the proper attention. To handle the transient and amplified oscillations in the UPD scheme, the non-uniform GPD schemes of both types would be a better choice to model the multiscale problems and to preserve the computational complexity of the method. Lee et al. \cite{2013_lee_variable} provided additional information regarding the advantages of employing variable mesoscopic time steps.

Unlike UPD scheme, in GPD scheme, one has the freedom to choose uniform or non-uniform distribution of GTTs and the extrapolation times based on the physical problem in time $\Delta t$. It is another significant advantage of the GPD schemes. 

In figures- \ref{fig:Contour_T1_NonUni} \& \ref{fig:Contour_T1_NonUni_Error}, GPD scheme of type-I is considered with non-uniform distribution of 24 number of GTTs and non-uniform distribution of extrapolation times within the long time step $\Delta t$. For example, a non-uniform distribution of 24 number of GTTs is considered as $k=\{4, 2, 8, 1, 5, 4\}$, and a non-uniform distribution of the extrapolation time is considered as $0.00076\times\{0.1, 0.15, 0.05, 0.25, 0.05, 0.4\}$, where $\Delta t-24\tau=0.00076$. From the contour plots of solutions and surface plot of errors, it is clear that GPD scheme of type-I is capable of dealing with variable GTTs and variable extrapolation times.

\begin{minipage}{0.39\linewidth}
	\centering
	\includegraphics[width=4.5cm,height=2.5cm]{AS_Contour_New.eps}\\
	\includegraphics[width=4.5cm,height=2.5cm]{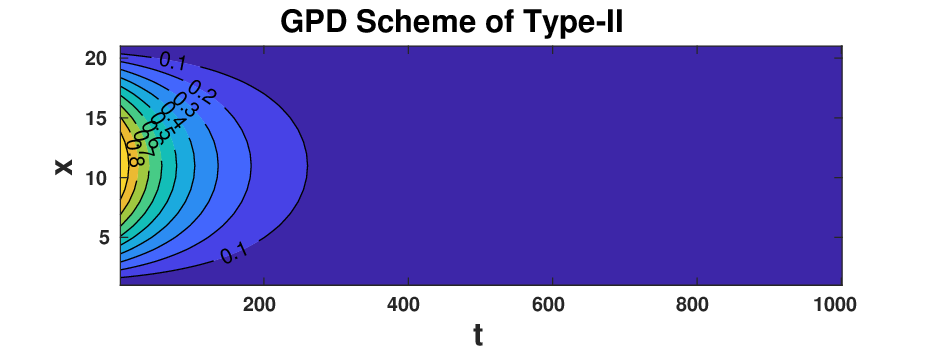}
	\captionof{figure}{Contours of the analytical solution and GPD solution of type-II for problem \eqref{eqn:22} with non-uniform distribution of GTTs and extrapolation times.}
	\label{fig:Contour_T2_NonUni}
\end{minipage}
\begin{minipage}{0.49\linewidth}
	\centering
	\includegraphics[width=5cm,height=4cm]{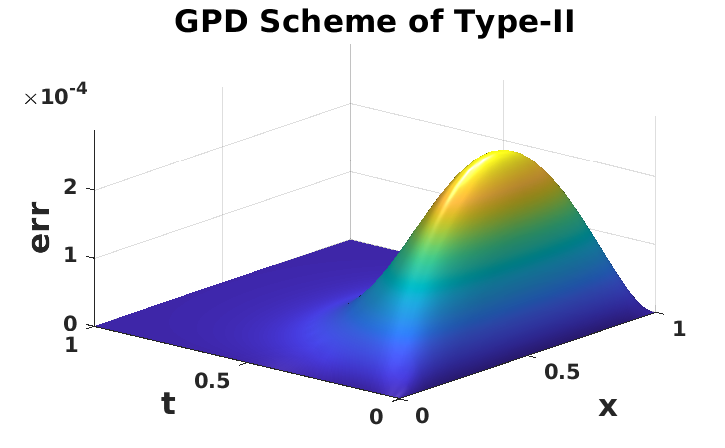}
	\captionof{figure}{Surface plots of errors for problem \eqref{eqn:22}.}
	\label{fig:Contour_T2_NonUni_Error}
\end{minipage}

In figures- \ref{fig:Contour_T2_NonUni} \& \ref{fig:Contour_T2_NonUni_Error}, GPD scheme of type-II is considered with non-uniform distribution of 24 number of GTTs and non-uniform distribution of extrapolation times within the time step $\Delta t$. For example, a non-uniform distribution of 24 number of GTTs is considered as $k=\{2, 7, 3, 2, 1, 4, 5\}$, and a non-uniform distribution of the extrapolation time is considered as $0.00076\times\{ 0.3, 0.05, 0.15, 0.25, 0.05, 0.2\}$, where $\Delta t-24\tau=0.00076$. From the contour plots of solutions and surface plot of error, it is clear that the GPD scheme of type-II is also capable of dealing efficiently with variable number of GTTs and variable extrapolation times. 

\section{Conclusion}

In this study, a generalized patch dynamics (GPD) scheme for multiscale problems is proposed. Based on the distribution of gap-tooth timesteppers (GTTs), GPD scheme has two different types, namely GPD scheme of type-I and type-II. GPD scheme of type-I is a `GTT - extrapolation' based modelling and GPD scheme of type-II is a `GTT - extrapolation - GTT' based modelling. Using these schemes, one can approximate an unavailable effective equation over macroscopic length and time scales, when only microscopic evolution laws are available. A microscopic model is only simulated over small patches of space-time domain using appropriate lifting operator and patch boundary conditions. In usual patch dynamics (UPD) scheme, the equal number of GTTs are implemented at the starting of each long time steps $\Delta t$ and consecutively one extrapolation is done in rest of the time step $\Delta t$. But in reality, based on the nature of the problem, one may need the flexibility to use less or more micro simulations at the high gradient zones in time and short or long extrapolations at low gradient zones in time \cite{2013_lee_variable}. This facilities are available in GPD scheme, but not in UPD scheme. Finer distribution of GTTs improves the solution compared to the coarser distribution of GTTs in GPD scheme. As shown in the model problem, the scheme is capable of providing good approximations for reaction-diffusion systems.

From this study, the following conclusions can be drawn

\begin{enumerate}
	\item GPD scheme is a generalized version of `UPD scheme', `Heterogeneous Multiscale Methods' (HMM), `Flow Averaging Integrators' (FLAVORS), `Variable Step Size Heterogeneous Multiscale Methods' (VSHMM), `Boosting Algorithm in Heterogeneous Multiscale Methods' (BA in HMM).
	\item The UPD scheme only uses two different time scales, namely micro and macro time scales. Whereas, the GPD scheme uses three different time scales: micro, meso and macro time scales to predict the system level behaviours.
	\item GPD scheme of both types are capable to give better accuracy and take less computational time compared to the UPD scheme (Tables \ref{table:1}, \ref{table:2}, \ref{table:3} and \ref{table:4}). The most effective approach to achieving better accuracy with less computational time compared to the UPD scheme, one should adopt to taking fewer macro time steps with a finer distribution of GTTs (i.e., more meso time steps) while still adhering to the relaxation time limit. 
	\item GPD scheme is more appropriately able to capture the physics of the problem compared to the UPD scheme. Non-uniform (variable) distributions of the GTTs as well as the extrapolation times based on the corresponding physical problem gives such additional features.
	
	\item GPD scheme of type-II is a modified version of GPD of type-I. In the GPD scheme of type-II, final microsimulations are executed just before termination of the macroscopic time step $\Delta t$. If the physical problem exhibits high gradient or oscillations near the terminating time in $\Delta t$ and requires microsimulations just before the macro time step $\Delta t$, then GPD-II will be a more approprite model than the UPD scheme with respect to the physical problem.
	
	\item Extrapolation time in UPD scheme has an upper bound, beyond that value UPD scheme fails to converge to the analytical solution, and the value of this upper bound depends on the corresponding model problems. Most of the times, GPD schemes are capable of overcoming such failures of the UPD scheme.
\end{enumerate}

\section*{Acknowledgments}
The first author thanks the Indian Institute of Technology Guwahati, India, for research space and infrastructure. Thanks to the anonymous reviewers for their constructive comments.

\bibliographystyle{siam}
\bibliography{references}
\end{document}